\documentclass[12 pt]{article}
\usepackage{graphicx} 
\usepackage[all]{xy}
\usepackage{pdflscape}
\usepackage[hypcap=false]{caption}
\usepackage{amsmath,amssymb,amsfonts,latexsym,float,graphics,epsfig} 
\usepackage{setspace}
\usepackage{xcolor}
\usepackage[top=2cm, bottom=2cm, left=2cm, right=2cm]{geometry}
\usepackage[colorlinks]{hyperref}
\usepackage{tikz}
\usetikzlibrary{cd}  

\newtheorem{theorem}{Theorem}[section]  

\newtheorem{proposition}[theorem]{Proposition}

\title{Locally Conformal Analysis of Non-Degenerate Odd Pairs}
\author{}
\date{}

\begin{document}
\maketitle

\begin{center}
Beg\"{u}m Ate\c{s}li$^{\ast,\ast\ast}$\footnote{e-mails: 
\href{mailto:begum.atesli@medipol.edu.tr}{begum.atesli@medipol.edu.tr}, \href{mailto:begumatesli@itu.edu.tr}{begumatesli@itu.edu.tr}, corresponding author}, Aybike \c{C}atal-\"{O}zer$^{\ast\ast}$\footnote{e-mail: 
\href{mailto:ozerayb@itu.edu.tr}{ozerayb@itu.edu.tr}}, and
O\u{g}ul Esen$^{\dagger,\ddagger}$\footnote{e-mail: \href{mailto:oesen@gtu.edu.tr}{oesen@gtu.edu.tr}}\\

  \bigskip

  $^{\ast}$School of Business, \\ İstanbul Medipol University, 34815 Kavacık-İstanbul, Turkey

\bigskip
  
$^{\ast\ast}$Department of Mathematics Engineering, \\ \.{I}stanbul Technical University,  34467 Maslak-\.{I}stanbul, Turkey \\

 \bigskip

  $^{\dagger}$Department of Mathematics, Gebze Technical University,  \\  41400 Gebze-Kocaeli, Turkey \\  

 \bigskip

  $^{\ddagger}$Center for Mathematics and its Applications,  \\ 
Khazar University, Baku, AZ1096, Azerbaijan
\end{center}

\bigskip

 \begin{abstract}
We study the locally conformal geometry of cosymplectic and coorientable contact manifolds within a common algebraic framework based on a non-degenerate odd pair: a one-form and a two-form whose appropriate wedge product defines a nowhere-vanishing top-degree form on an odd-dimensional manifold. We consider conformally compatible local geometric and Hamiltonian data and encode their gluing through a closed Lee form and the Lichnerowicz–de Rham differential.

For locally conformal cosymplectic geometry, we introduce an equal-weight formulation in which the defining forms carry the same conformal weight. This formulation allows compatible local Hamiltonian vector fields to glue directly into a global Hamiltonian vector field without additional rescaling, thereby preserving the local dynamics. For contact geometry, we develop a twisted odd-pair formulation in which the Lee form records the conformal gluing of local Hamiltonian data, while the global one-form remains an ordinary contact form. We illustrate these constructions through time-dependent mass dynamics in the cosymplectic setting and dissipative particle dynamics in the contact setting.

\end{abstract}

\smallskip

 
\tableofcontents
\onehalfspacing

\numberwithin{equation}{section}

\setlength{\parindent}{2em}
\setlength{\parskip}{3ex}

\section{Introduction}

Poisson geometry provides a natural framework for Hamiltonian dynamics \cite{Laurent13,Vaisman94,Weinstein98}, with applications throughout classical mechanics and mathematical physics \cite{AbrahamMarsden,Arn,LiMa,MaRa99}. The Jacobi identity ensures that Hamiltonian vector fields are closed under the Jacobi--Lie bracket, while the skew-symmetry of the Poisson bracket implies that each Hamiltonian function is conserved along its own Hamiltonian flow.

Geometric descriptions of nonconservative systems motivate extensions of the Poisson framework. One such extension is Jacobi geometry \cite{Lichnerowicz-Jacobi,Marle91,vaisman2002jacobi}, originating in Kirillov's notion of a local Lie algebra \cite{Kirillov}. Its bracket retains the Jacobi identity but satisfies a first-order differential property in place of the Poisson Leibniz rule. In this work, we approach the relation between Poisson and Jacobi geometries through differential forms. Starting from a common algebraic framework, we examine how conformal modifications of the defining differential relations affect the associated brackets and Hamiltonian dynamics. For related approaches based on differential forms and pairs, see \cite{de2026poisson,gracia2026pairs}.

\noindent
\textbf{Non-degenerate Odd Pairs.}
The common algebraic starting point of our analysis is a non-degenerate odd
pair on a \((2n+1)\)-dimensional manifold. Such a pair consists of a
one-form \(\kappa\) and a two-form \(\Omega\) satisfying
\[
\kappa\wedge\Omega^n\neq0.
\]
This algebraic condition determines a distinguished Reeb vector field, a
bivector field, and the corresponding musical mappings. Cosymplectic and
coorientable contact geometries arise as two different differential
realizations of this common datum. Throughout this paper, all contact
manifolds are assumed to be coorientable; hence, for brevity, we subsequently
use the term ``contact manifold'' without the qualifier ``coorientable.''

A cosymplectic pair \((\tau,\Omega)\), characterized by
\(d\tau=0\) and \(d\Omega=0\), induces a Poisson structure. By contrast, a
contact form \(\eta\), together with the two-form \(\Omega=d\eta\), induces a
Jacobi structure. The odd-pair formulation therefore allows the Poisson
character of cosymplectic geometry and the Jacobi character of contact
geometry to be studied within a unified framework.

\noindent \textbf{Locally Conformal Geometry.}
Our main concern is the passage from local to global geometry. While globally defined tensorial structures are directly reflected in local charts, the converse is more subtle: locally defined tensorial quantities do not necessarily determine a global tensorial object, unless appropriate compatibility conditions are satisfied on chart overlaps.
More concretely, the existence
of cosymplectic or contact data on every chart does not by itself imply the
existence of a corresponding global structure.

In the present work, we focus
on situations in which the local forms fail to glue in the ordinary
tensorial sense but become compatible after conformal rescaling. The
derivatives of the conformal factors then determine a closed Lee form
\(\varphi\), and the global differential relations are naturally expressed
through the Lichnerowicz--de~Rham differential
\cite{GuLi84},
\[
d_\varphi=d-\varphi\wedge.
\]
Thus, the twisted differential is not introduced merely as a formal
replacement for the exterior derivative; it records the obstruction to the
ordinary globalization of the local geometric data.

The symplectic prototype of this construction is the theory of locally
conformally symplectic manifolds, initiated by Vaisman \cite{Vaisman85} and
subsequently developed in several geometric and dynamical directions; see,
for example,
\cite{Banyaga,Bazzoni2018,EsenLeonSarZaj1,esen2023conformal,
EsLaLeMa21,HallerRybicki,OtSt15}. Although the locally conformally
symplectic case is not treated separately here, it provides the guiding
principle: a collection of locally closed structures becomes globally
meaningful after its conformal mismatch is encoded by a Lee form. We carry
this philosophy to the odd-dimensional cosymplectic and coorientable contact
settings.

\noindent
\textbf{Goal I: An Equal-Weight Formulation of Locally Conformally Cosymplectic Geometry.}
Cosymplectic manifolds provide a geometric framework for Hamiltonian systems with explicit time dependence \cite{Cape,leonsardon1,LeonTuyn}. Their locally conformal counterparts have been studied, in particular, in \cite{AtEsLe23,ChLeMa91}, where the defining one-form and two-form carry different conformal weights.

Here we consider an equal-weight formulation in which both forms are governed by the same closed Lee form and the same Lichnerowicz--de~Rham differential. We establish its local-to-global characterization and determine the associated Jacobi structure. When the local Hamiltonians carry the same conformal weight as the defining two-form, the local Hamiltonian vector fields glue directly into a global vector field. This construction agrees with the dynamics obtained by replacing the ordinary differential of the Hamiltonian with its Lichnerowicz--de~Rham differential.

We treat Hamiltonian and evolution dynamics separately. The Hamiltonian construction applies to a general closed Lee form, whereas the gluing of locally normalized evolution fields requires an additional compatibility condition. In the scalar formulation considered here, this condition is expressed through a nowhere-vanishing global conformal unit and restricts the evolution construction to the exact-Lee sector. We also establish the closure of the Hamiltonian vector fields under the Jacobi--Lie bracket and the integrability of their characteristic distribution. These results concern the geometric integrability of the distribution, rather than complete integrability in the Liouville sense. A particle with a time-dependent mass tensor illustrates the resulting equations.

\noindent
\textbf{Goal II: A Conformally Twisted Odd-Pair Formulation of Contact Geometry.}
Contact geometry provides a complementary framework for Hamiltonian descriptions of dissipative systems and thermodynamics \cite{arnold98,Bravetti2017,Bravettithermo,de2019contact,EsGu11,Grmela14}. A contact structure has an intrinsic conformal freedom: multiplication of a contact form by a nowhere-vanishing function preserves its kernel. We examine how this freedom enters the globalization of local contact Hamiltonian systems.

For a global one-form \(\eta\) and a closed Lee form \(\varphi\), we consider the twisted odd pair determined by \(\Omega=d_\varphi\eta\). Its non-degeneracy is equivalent to the ordinary contact condition on \(\eta\). Moreover, we show that the associated bivector and the Lee-corrected odd-pair Reeb vector field recover the canonical Jacobi structure of \(\eta\). The construction therefore provides a twisted odd-pair description of ordinary coorientable contact geometry.

Within this description, we prove that direct Lichnerowicz--de~Rham deformation and conformal globalization of compatible local contact Hamiltonian systems yield the same global dynamics. We derive the corresponding coordinate equations and illustrate them through a dissipative particle model, making explicit the Lee-form-dependent terms arising in the chosen representation.

\noindent
\textbf{Organization of the Paper.}
Section~\ref{LC-Sec-1} recalls the Lichnerowicz--de~Rham differential and introduces non-degenerate odd pairs. Section~\ref{LC-Sec-121} reviews the Poisson structure associated with cosymplectic geometry and develops its equal-weight locally conformal formulation, including Hamiltonian dynamics, evolution dynamics in the exact-Lee sector, and the time-dependent mass example. Section~\ref{LC-Sec-12} develops the twisted odd-pair formulation of contact geometry, establishes its relation to the canonical contact Jacobi structure, and studies the globalization of local Hamiltonian dynamics and the dissipative particle example. Section~\ref{sec:summary-discussion} discusses the results and possible extensions. The appendix collects the Schouten--Nijenhuis and right-contraction identities used in the paper.

\section{Fundamentals}\label{LC-Sec-1}

This section introduces the two ingredients used in the subsequent constructions. The Lichnerowicz--de~Rham differential encodes the conformal compatibility of local differential forms, while a non-degenerate odd pair provides the common algebraic framework for cosymplectic and contact geometries. We also clarify the additional condition required for a global scalar conformal unit, which will enter the evolution dynamics.

\subsection{The Lichnerowicz--de~Rham Differential}\label{sec-LdR}

Let \(M\) be a smooth manifold with a suitable open cover
\(\{U_\alpha\}\). Suppose that each chart \(U_\alpha\) is equipped with a
smooth function \(\sigma_\alpha\) such that their exterior derivatives agree
on every nonempty overlap:
\begin{equation}\label{Lee-potential-compatibility}
d\sigma_\alpha
=
d\sigma_\beta
\qquad
\text{on }
U_\alpha\cap U_\beta.
\end{equation}
Consequently, the local one-forms \(d\sigma_\alpha\) glue to a globally
defined one-form \(\varphi\), characterized locally by
\begin{equation}\label{local-Lee-potential}
\varphi|_\alpha=d\sigma_\alpha.
\end{equation}
Since
\begin{equation}
d\varphi|_\alpha=d^2\sigma_\alpha=0
\end{equation}
on every chart, \(\varphi\) is closed. This globally defined closed one-form
is called the Lee form. Moreover,
\eqref{Lee-potential-compatibility} implies
\begin{equation}
\sigma_\alpha-\sigma_\beta=c_{\alpha\beta}
\qquad
\text{on }
U_\alpha\cap U_\beta,
\end{equation}
where \(c_{\alpha\beta}\) is locally constant.

The Lichnerowicz--de~Rham differential associated with the Lee form
\(\varphi\) \cite{GuLi84} is obtained by deforming the ordinary exterior
derivative through the wedge action of \(\varphi\). More precisely, it is
the operator
\begin{equation}\label{LdR-Diff-}
d_\varphi:\Gamma^k(M)\longrightarrow\Gamma^{k+1}(M),
\qquad
\beta\longmapsto d\beta-\varphi\wedge\beta.
\end{equation}
For any differential form \(\beta\), a direct computation gives
\begin{equation}
d_\varphi^2\beta=-d\varphi\wedge\beta.
\end{equation}
Thus, the closedness of \(\varphi\) implies \(d_\varphi^2=0\).
The cohomology of this complex is called the Lichnerowicz--de~Rham
cohomology \cite{HaRy99}.

The relation with the ordinary exterior derivative is particularly
transparent on a chart where \(\varphi|_\alpha=d\sigma_\alpha\):
\begin{equation}\label{LdR-local-conjugation}
\left.d_\varphi\beta\right|_\alpha
=
e^{\sigma_\alpha}
d\left(e^{-\sigma_\alpha}\beta|_\alpha\right).
\end{equation}
Consequently, a global form is \(d_\varphi\)-closed precisely when its
locally rescaled representatives are closed. This identity will underlie
the local-to-global constructions below.

\noindent
\textbf{The Scalar Conformal Unit.}
For the ordinary constant function \(1\), the twisted differential gives
\begin{equation}
d_\varphi 1=-\varphi.
\end{equation}
Thus, the usual unit is \(d_\varphi\)-closed only when the Lee form vanishes.
In the evolution dynamics considered below, we instead use a nowhere-vanishing
function \(\mathbb I\) satisfying
\begin{equation}\label{twisted-closed-I}
d_\varphi\mathbb I=0.
\end{equation}
We call such a function a scalar conformal unit. Its existence is an
additional global condition: it holds if and only if the Lee form is exact.

Indeed, a nowhere-vanishing solution of \eqref{twisted-closed-I} satisfies
\begin{equation}
\varphi
=
\frac{d\mathbb I}{\mathbb I}
=
d\log|\mathbb I|.
\end{equation}
Conversely, if \(\varphi=d\sigma\) globally, then
\(\mathbb I=e^\sigma\) is a nowhere-vanishing solution.

The local description also makes the gluing condition explicit. On a chart
with \(\varphi|_\alpha=d\sigma_\alpha\), every nowhere-vanishing local
solution has the form
\begin{equation}\label{local-I}
\mathbb I|_\alpha=e^{\sigma_\alpha}1_\alpha,
\end{equation}
where \(1_\alpha\) is a nonzero locally constant function. Here
\(1_\alpha\) denotes a local normalization, which need not equal the
numerical value \(1\). These representatives define a global scalar
conformal unit precisely when
\begin{equation}\label{local-unit-gluing}
e^{\sigma_\alpha}1_\alpha
=
e^{\sigma_\beta}1_\beta
\qquad
\text{on }U_\alpha\cap U_\beta.
\end{equation}
Their twisted closedness follows directly from
\begin{equation}\label{local-twisted-closed-I}
\left.d_\varphi\mathbb I\right|_\alpha
=
e^{\sigma_\alpha}d1_\alpha
=
0.
\end{equation}

For a non-exact Lee form, such normalizations can be chosen locally but
cannot satisfy \eqref{local-unit-gluing} globally. Accordingly, the scalar
conformal unit will be used only in the exact-Lee sector.

\subsection{Non-degenerate Odd Pair}\label{sec:ndop}

Let \(M\) be a \((2n+1)\)-dimensional manifold. A pair
\((\kappa,\Omega)\), consisting of a one-form \(\kappa\) and a two-form
\(\Omega\), is called a non-degenerate odd pair if
\begin{equation}\label{non-deg-odd-pair}
    \kappa\wedge\Omega^n\neq 0.
\end{equation}
This condition implies that \(\Omega\) has rank \(2n\) and that its one-dimensional kernel is transverse to \(\ker\kappa\). Equivalently, the restriction of \(\Omega\) to \(\ker\kappa\) is non-degenerate. These algebraic properties determine the musical mappings, Reeb vector field, and bivector field introduced below, without imposing any differential conditions on \(\kappa\) or \(\Omega\).

This non-degeneracy condition \eqref{non-deg-odd-pair} induces a musical isomorphism $\flat^{\kappa}$ from the space $\mathfrak{X}(M)$ of vector fields to the space $\Gamma^1(M)$ of one-form sections:
\begin{equation}\label{flat-map-kappa}
\flat^{\kappa}:\mathfrak{X}(M)\longrightarrow \Gamma^1(M),
\qquad
X\longmapsto \iota_X\Omega+(\iota_X\kappa)\kappa.
\end{equation}
To see that this map is invertible, suppose that
\(\flat^\kappa(X)=0\). Evaluation on a nonzero vector in \(\ker\Omega\) gives \(\kappa(X)=0\), and hence \(\iota_X\Omega=0\). Since \(\ker\Omega\cap\ker\kappa={0}\), it follows that \(X=0\). The map is therefore an isomorphism at each point.
Since $\flat^{\kappa}$ is an isomorphism, it is invertible; we denote its inverse by
\begin{equation}\label{sharp-kappa}
\sharp^{\kappa}: \Gamma^1(M) \longrightarrow \mathfrak{X}(M),
\qquad
\sharp^{\kappa}=(\flat^{\kappa})^{-1}.
\end{equation}
In this geometry, there exists a distinguished vector field, called the Reeb
vector field \(\mathcal{R}^{\kappa}\), defined by
\begin{equation}\label{Reeb-kappa}
\iota_{\mathcal{R}^{\kappa}}\kappa=1,\qquad
\iota_{\mathcal{R}^{\kappa}}\Omega=0,
\qquad \text{or equivalently,} \qquad
\mathcal{R}^{\kappa}=\sharp^{\kappa}\kappa.
\end{equation}
For one-forms \(\alpha,\beta\), we define the associated bivector field \(\Lambda^\kappa\) by
\begin{equation}\label{Lambda-kappa}
\Lambda^{\kappa} (\alpha,\beta)
=
\Omega(\sharp^\kappa \alpha, \sharp^\kappa \beta).
\end{equation}
The associated sharp mapping of \(\Lambda^{\kappa}\), obtained by right
contraction, is given by
\begin{equation}\label{sharp-lambda-kappa}
\sharp_{\Lambda^{\kappa}}:\Gamma^1(M)\longrightarrow \mathfrak{X}(M),
\qquad
\beta \longmapsto
\iota_{\beta}\Lambda^{\kappa}
=
\Lambda^{\kappa}(\bullet,\beta)
=
\sharp^{\kappa}\beta
-
\langle \beta,\mathcal{R}^{\kappa}\rangle \mathcal{R}^{\kappa}.
\end{equation}
This mapping is not an isomorphism; its kernel is generated by \(\kappa\).
  
To summarize, starting from the non-degenerate odd pair
\((\kappa,\Omega)\), we have canonically obtained the Reeb vector field
\(\mathcal{R}^{\kappa}\) and the bivector field \(\Lambda^{\kappa}\). Thus, the covariant
data \((\kappa,\Omega)\) determine the contravariant data
\((\Lambda^{\kappa},\mathcal{R}^{\kappa})\). Conversely, if \(\Lambda^{\kappa}\) has rank \(2n\)
and \(\mathcal R^{\kappa}\) is transverse to the characteristic distribution
\(\operatorname{Im}\sharp_{\Lambda^{\kappa}}\), then \((\Lambda^{\kappa},\mathcal R^{\kappa})\) determines
\((\kappa,\Omega)\) uniquely by requiring
\begin{equation}\label{dual-data-reconstruction}
\kappa(\mathcal R ^{\kappa})=1,
\qquad
\ker\kappa=\operatorname{Im}\sharp_{\Lambda^{\kappa}},
\qquad
\iota_{\mathcal R ^{\kappa}}\Omega=0,
\end{equation}
and by taking \(\Omega\) to be the inverse of \(\Lambda^{\kappa}\) on
\(\operatorname{Im}\sharp_{\Lambda^{\kappa}}\). In this sense, under the appropriate
non-degeneracy condition, the covariant data \((\kappa,\Omega)\) and the
contravariant data \((\Lambda ^{\kappa},\mathcal R ^{\kappa})\) are equivalent descriptions of the
same odd-dimensional geometry. 

The subsequent sections develop two bracket constructions from this common algebraic datum. The bivector \(\Lambda^\kappa\) alone defines an almost Poisson bracket, while the pair \((\Lambda^\kappa,\mathcal R^\kappa)\) defines an almost Jacobi bracket. Neither construction satisfies the Jacobi identity in general. We examine the differential conditions that yield Poisson structures in the cosymplectic case and Jacobi structures in the contact case, and then study their locally conformal formulations.

\section{Locally Conformal Analysis of Cosymplectic Framework}\label{LC-Sec-121} 

In this section, we study cosymplectic manifolds in light of the non-degenerate odd-pair framework introduced in Subsection~\ref{sec:ndop}. Then we present the locally conformal analysis of the cosymplectic geometry. Starting from the ordinary cosymplectic setting recalled above, we replace the de~Rham closedness conditions on \(\tau\) and \(\Omega\) by their Lichnerowicz--de~Rham analogues. This leads to an equal-weight locally conformally cosymplectic structure, in which both forms are controlled by the same Lee one-form. We then examine the local-to-global meaning of this definition, derive the corresponding Jacobi structure, and formulate the Hamiltonian and evolution dynamics associated with the locally conformal cosymplectic data.
 
\subsection{Poisson Character of an Odd Pair: Cosymplectic Manifolds}
\label{sec:cos}

Let \(M\) be a \((2n+1)\)-dimensional manifold equipped with a
non-degenerate odd pair \((\tau,\Omega)\), as defined in
Subsection~\ref{sec:ndop}. Thus, \(\tau\) is a one-form and \(\Omega\) is a
two-form satisfying the nondegeneracy condition
\begin{equation}
\tau\wedge\Omega^n\neq0.
\end{equation}
According to the discussion in Subsection~\ref{sec:ndop}, the odd pair \((\tau,\Omega)\) may also be described by the corresponding contravariant
data \((\Lambda^\tau,\mathcal R^\tau)\). Starting from \((\tau,\Omega)\), the Reeb
vector field \(\mathcal R^\tau\) is determined by the conditions
\begin{equation}\label{Reeb-tau}
\iota_{\mathcal R^\tau}\tau=1,
\qquad
\iota_{\mathcal R^\tau}\Omega=0,
\end{equation}
whereas the bivector field \(\Lambda^\tau\) is defined by
\begin{equation}\label{Lambda-tau}
\Lambda^\tau(dF,dH)
=
\Omega(\sharp^\tau dF,\sharp^\tau dH).
\end{equation}
Conversely, under the corresponding rank and transversality assumptions, the
pair \((\Lambda^\tau,\mathcal R^\tau)\) recovers the original covariant data
\((\tau,\Omega)\). Thus, in what follows, we freely pass between the
covariant description \((\tau,\Omega)\) and the contravariant description
\((\Lambda^\tau,\mathcal R^\tau)\).

We now define a skew-symmetric bracket on \(C^\infty(M)\) by using only the
bivector field \(\Lambda^\tau\). More precisely, for \(F,H\in C^\infty(M)\), we set
\begin{equation}\label{almost-poisson-bracket-tau}
    \{F,H\}
    =
    \Lambda^\tau(dF,dH).
\end{equation}
This bracket is almost Poisson in the sense that it is skew-symmetric and
satisfies the Leibniz identity: 
\begin{equation}\label{almost-poisson-leibniz-tau}
\{F,KH\}
=
K\{F,H\}
+
H\{F,K\}.
\end{equation}
However, for a general
non-degenerate odd pair, the Jacobi identity is not imposed. 
Thus, the bracket behaves as a derivation with respect to the pointwise
multiplication of functions, but it does not necessarily define a Lie algebra
structure on \(C^\infty(M)\).

The corresponding almost Poisson Hamiltonian vector field is obtained by
contracting \(dH\) with the bivector field \(\Lambda^\tau\):
\begin{equation}\label{almost-poisson-Ham-vf-tau}
X_H
=
\sharp_{\Lambda^\tau}(dH)
=
\sharp^\tau(dH)
-
\mathcal R^\tau(H)\mathcal R^\tau .
\end{equation}
Equivalently, \(X_H\) may be characterized intrinsically by
\begin{equation}\label{almost-poisson-Ham-vf-char-tau}
\iota_{X_H}\Omega
=
dH-\mathcal R^\tau(H)\tau,
\qquad
\iota_{X_H}\tau
=
0.
\end{equation}
The second condition shows that \(X_H\) is tangent to the distribution
\(\ker\tau\). 

 The failure of the almost Poisson bracket to satisfy the Jacobi identity has
a direct dynamical consequence: the associated Hamiltonian vector fields need
not be closed under the Jacobi--Lie bracket.

\begin{proposition}\label{prop:almost-poisson-Ham-closure}
Let \((M,\tau,\Omega)\) be a non-degenerate odd pair, let
\(\Lambda^\tau\) be its associated bivector field, and denote the
Hamiltonian vector field of \(H\in C^\infty(M)\) by \(X_H\). Then, for every
\(F,H,K\in C^\infty(M)\),
\begin{equation}\label{almost-poisson-Ham-vf-bracket-defect-tau}
[X_F,X_H](K)
=
-X_{\{F,H\}}(K)
+
\frac{1}{2}
[\Lambda^\tau,\Lambda^\tau]_{SN}(dF,dH,dK).
\end{equation}
\end{proposition}

Accordingly, the Schouten--Nijenhuis bracket
\([\Lambda^\tau,\Lambda^\tau]_{SN}\) measures the extent to which the Hamiltonian assignment fails to be a Lie algebra anti-homomorphism. When this
tensor does not vanish, the family of Hamiltonian vector fields may fail to
close under the Jacobi--Lie bracket. We refer the reader to Appendix \ref{sec:SN} for some details of the Schouten--Nijenhuis bracket.  

The Hamiltonian function is conserved along its own
Hamiltonian flow:
\begin{equation}\label{Hamiltonian-conservation-XH-tau}
X_H(H)=0.
\end{equation} Let \(\mathrm{Vol}=\tau\wedge\Omega^n\) be the volume form determined by the
non-degenerate odd pair. Although \(H\) is conserved along \(X_H\), the
corresponding volume form need not be preserved. Indeed, by Cartan's formula,
we compute
\begin{equation}\label{volume-XH-tau}
\begin{split}
\mathcal L_{X_H}\mathrm{Vol}
={}&
(\iota_{X_H}d\tau)\wedge\Omega^n
+
n\tau\wedge(\iota_{X_H}d\Omega)\wedge\Omega^{n-1}
-
n\mathcal R^\tau(H)\tau\wedge d\tau\wedge\Omega^{n-1}.
\end{split}
\end{equation}
Therefore, \(\mathrm{Vol}\) is not preserved in general. Its preservation is
an additional geometric property, rather than a consequence of the
non-degeneracy condition alone.

\noindent 
\textbf{Cosymplectic Manifolds.}
Suppose that the differential forms determining the
non-degenerate odd pair are closed, namely
\begin{equation}\label{closure-cosym-tau}
\tau\wedge\Omega^n\neq0,
\qquad
d\tau=0,
\qquad
d\Omega=0.
\end{equation}
Then \((M,\tau,\Omega)\) is a cosymplectic manifold \cite{Cape}. The
closedness conditions eliminate the Schouten--Nijenhuis obstruction
identified above. This observation leads to the following result \cite{Cape}.

\begin{proposition}\label{prop:cosymplectic-Poisson}
Let \(M\) be a \((2n+1)\)-dimensional manifold equipped with a
non-degenerate odd pair \((\tau,\Omega)\) satisfying closedness conditions 
\eqref{closure-cosym-tau}. Then the associated bivector field
\(\Lambda^\tau\) is Poisson:
\begin{equation}\label{Poisson-ident-tau}
[\Lambda^\tau,\Lambda^\tau]_{SN}=0.
\end{equation}
Consequently, the bracket \eqref{almost-poisson-bracket-tau} 
satisfies the Jacobi identity
\begin{equation}\label{Jac-id-Poisson-tau}
\circlearrowright_{F,H,K}
\{F,\{H,K\}\}
=
0
\end{equation}
and therefore defines a Poisson bracket on \(C^\infty(M)\).
\end{proposition}

It now follows from
\eqref{almost-poisson-Ham-vf-bracket-defect-tau} that, when the bivector field
satisfies the Poisson condition \eqref{Poisson-ident-tau}, namely,
\(
[\Lambda^\tau,\Lambda^\tau]_{SN}=0,
\)
the defect term on the right-hand side of
\eqref{almost-poisson-Ham-vf-bracket-defect-tau} vanishes identically.
Hence, the cosymplectic Hamiltonian vector fields satisfy
\begin{equation}\label{Poisson-integ-tau}
[X_F,X_H]
=
-X_{\{F,H\}}.
\end{equation}
Therefore, the space
\[
\mathfrak X_{\mathrm{Ham}}(M)
=
\left\{X_H: H\in C^\infty(M)\right\}
\]
is closed under the Jacobi--Lie bracket and hence forms a Lie algebra.

Furthermore, the characteristic distribution of the Poisson bivector,
\begin{equation}
\mathcal D_{\Lambda^\tau}
=
\operatorname{Im}\sharp_{\Lambda^\tau}
=
\operatorname{span}
\left\{
X_H:\ H\in C^\infty(M)
\right\}
=
\ker\tau,
\end{equation}
is involutive, since \eqref{Poisson-integ-tau} shows that the Jacobi--Lie
bracket of any two Hamiltonian vector fields is again Hamiltonian.
It is therefore integrable. The integral leaves are precisely the
symplectic leaves of the Poisson structure \(\Lambda^\tau\), in accordance
with the Stefan--Sussmann theorem
\cite{BhaskaraViswanath,DuZu05,Laurent13,Vaisman94,Weinstein98}.

In the cosymplectic case, we still have the conservation of the Hamiltonian
function along its Hamiltonian vector field, namely,
$X_H(H)=0$. Moreover, a stronger preservation property holds in the cosymplectic case:
the volume form is also preserved by the Hamiltonian flow.
Indeed, by \eqref{volume-XH-tau}, we have
\begin{equation}\label{cosymp-volume-preservation-tau}
\mathcal L_{X_H}\mathrm{Vol}=0.
\end{equation}

\noindent \textbf{Coordinate Formulas.} 
The cotangent bundle \(T^*Q\) admits the exact symplectic structure
\begin{equation}
    \Omega_Q=-d\Theta_Q,
\end{equation}
where \(\Theta_Q\) is the canonical one-form.
Accordingly, the extended cotangent bundle \(\mathbb R\times T^*Q\) can be
regarded as a cosymplectic manifold. Let \(s\) denote the global coordinate on
\(\mathbb R\). Pulling back the canonical symplectic two-form \(\Omega_Q\) on
\(T^*Q\) to \(\mathbb R\times T^*Q\), we obtain the canonical cosymplectic
structure
\begin{equation}\label{cancosypair}
(\mathbb R\times T^*Q,\tau,\Omega),
\qquad
\tau=ds,
\qquad
\Omega=\Omega_Q=-d\Theta_Q.
\end{equation}
Here, in order to avoid notational inflation, we use the same symbols for the
canonical forms on \(T^*Q\) and for their pullbacks to
\(\mathbb R\times T^*Q\). 

On a chart \(U_\alpha\), let \((s,q^i,p_i)\) denote the Darboux coordinates. Then we have
\begin{equation}\label{cancosypair-coord}
\tau \vert_{\alpha}=ds,
\qquad
\Theta_Q \vert_{\alpha}=p_i\,dq^i,
\qquad
\Omega \vert_{\alpha}=dq^i\wedge dp_i.
\end{equation}
Accordingly, the Reeb vector field is given by
\begin{equation}\label{Reeb-canonical-cosymp}
\mathcal R^\tau \vert_{\alpha}
=
\frac{\partial}{\partial s}.
\end{equation}
The cosymplectic bivector field \(\Lambda^\tau\) takes the local form
\begin{equation}\label{AP-Bra}
\Lambda^\tau\vert_{\alpha}
=
\frac{\partial}{\partial q^{i}}\wedge \frac{\partial}{\partial p_{i}}
.
\end{equation} For a Hamiltonian function \(H=H(s,q,p)\), the Hamiltonian vector field
\(X_H\) is
\begin{equation}\label{Hvf-cosymp-local-vector}
X_H \vert_{\alpha}
=
\frac{\partial H}{\partial p_i}\frac{\partial}{\partial q^i}
-
\frac{\partial H}{\partial q^i}\frac{\partial}{\partial p_i}.
\end{equation}
Thus, along an integral curve of \(X_H\), the coordinate \(s\) is constant,
since \(X_H(s)=0\). Therefore, \(X_H\) describes the Hamiltonian dynamics on
each fixed-time slice \(\{s=\mathrm{constant}\}\). Equivalently, with respect
to an arbitrary curve parameter \(t\), one has
\begin{equation}\label{Hvf-cosymp-local}
\frac{ds}{dt}=0,
\qquad
\frac{dq^i}{dt}
=
\frac{\partial H}{\partial p_i}(s,q,p),
\qquad
\frac{dp_i}{dt}
=
-\frac{\partial H}{\partial q^i}(s,q,p).
\end{equation}

\noindent \textbf{Cosymplectic Evolution Dynamics.} In order to formulate the corresponding non-autonomous Hamiltonian system, we
shall regard the additional parameter \(s\) as the flow parameter of the
dynamics. Consequently, \(s\) will play the role of the time variable. Accordingly, one introduces the evolution vector field
\(\mathcal E_H\), defined by
\begin{equation}\label{evo-cosymp}
\iota_{\mathcal E_H}\tau=1,
\qquad
\iota_{\mathcal E_H}\Omega
=
dH-\mathcal R^\tau(H)\tau.
\end{equation}

Equivalently, the evolution vector field is obtained by adding the Reeb vector
field to the Hamiltonian vector field. In other words, while \(X_H\) describes
the Hamiltonian dynamics on the fixed-\(s\) slices, the Reeb contribution
\(\mathcal R^\tau\) accounts for the evolution of the parameter \(s\) itself.
Thus,
\begin{equation}\label{evo-cosymp-vector}
\mathcal E_H
=
\mathcal R^\tau+X_H.
\end{equation}
Note that the conservation of \(H\) along the Hamiltonian vector field
\(X_H\) does not imply its conservation along the evolution vector field.
Indeed, one has, in general,
\begin{equation}
\mathcal E_H(H)=\mathcal R^\tau(H),
\end{equation}
so that \(H\) is conserved along the non-autonomous evolution only when it is
independent of the parameter \(s\).

In local Darboux coordinates \((s,q^i,p_i)\), the evolution  vector field is computed to be
\begin{equation}\label{evo-cosymp-vector-local}
\mathcal E_H \vert_{\alpha}
=
\frac{\partial}{\partial s}
+
\frac{\partial H}{\partial p_i}\frac{\partial}{\partial q^i}
-
\frac{\partial H}{\partial q^i}\frac{\partial}{\partial p_i}.
\end{equation}
Then the evolution equations take the following form:
\begin{equation}\label{evo-cosymp-local-lambda}
\frac{ds}{dt}=1,
\qquad
\frac{dq^i}{dt}
=
\frac{\partial H}{\partial p_i}(s,q,p),
\qquad
\frac{dp_i}{dt}
=
-\frac{\partial H}{\partial q^i}(s,q,p).
\end{equation}
Hence, \(s(t)=t+s_0\). After choosing the origin of the parameter
so that \(s=t\), the geometric coordinate \(s\) plays the role of the
physical time variable, and the evolution equations take the  
non-autonomous Hamiltonian form
\begin{equation}\label{evo-cosymp-local-time}
\frac{dq^i}{ds}
=
\frac{\partial H}{\partial p_i}(s,q,p),
\qquad
\frac{dp_i}{ds}
=
-\frac{\partial H}{\partial q^i}(s,q,p).
\end{equation}
Thus, in the cosymplectic description, \(X_H\) is the Hamiltonian vector field
on fixed-\(s\) slices, whereas \(\mathcal E_H=\mathcal R^\tau+X_H\) is the
evolution vector field whose integral curves may be parametrized by the time
coordinate \(s\).

\subsection{Conformal Analysis of Cosymplectic Manifolds}

Consider a non-degenerate odd-dimensional manifold \((M,\tau,\Omega)\), as in
Subsection~\ref{sec:cos}, together with the bivector \(\Lambda^\tau\) defined
in \eqref{Lambda-tau} and the Reeb vector field \(\mathcal R^\tau\) defined in
\eqref{Reeb-tau}. In Section~\ref{sec:cos}, we observed that the integrability
of the associated Hamiltonian vector fields is controlled by the de~Rham
closedness of \(\tau\) and \(\Omega\). More precisely, the cosymplectic conditions \(d\tau=0\) and \(d\Omega=0\) are sufficient to ensure that the associated bivector \(\Lambda^\tau\) defines a Poisson bracket.

We now relax the  closedness conditions  \(d\tau=0\) and \(d\Omega=0\) by replacing the
ordinary de~Rham differential with the Lichnerowicz--de~Rham differential.
This conformal relaxation naturally turns the Poisson bracket into a Jacobi
bracket: the resulting bracket still satisfies the Jacobi identity, but in general it is no longer a derivation in each argument with respect to the usual pointwise product.

Accordingly, we call a non-degenerate odd-dimensional manifold
\((M,\tau,\Omega)\) a locally conformally cosymplectic manifold, abbreviated as
an LCCos manifold, if there exists a closed one-form \(\varphi\), called the
Lee one-form, such that
\begin{equation}\label{LCC-str}
\tau\wedge\Omega^n\neq 0, \qquad
d_{\varphi}\tau=0,\quad
d_{\varphi}\Omega=0.
\end{equation}
Here \(d_{\varphi}\) denotes the Lichnerowicz--de~Rham differential. The use of
the term ``locally conformally'' will be justified below by showing that
\eqref{LCC-str} is equivalent, locally, to an ordinary cosymplectic structure
after a conformal rescaling. We denote an LCCos manifold by \((M,\tau,\Omega;\varphi)\).

Using the bivector field \(\Lambda^\tau\), the Lee one-form \(\varphi\) gives
rise to a distinguished vector field by contraction. More precisely, we define
\begin{equation}\label{Lee-vector-LCCos}
\mathcal{V}_{\varphi}
=
\iota_{\varphi}\Lambda^\tau
=
\Lambda^\tau(\bullet,\varphi).
\end{equation}
We shall refer to \(\mathcal{V}_{\varphi}\) as the Lee vector field associated
with the LCCos structure. It records the component of the Lee form detected by the bivector \(\Lambda^\tau\). Indeed, when
\(\varphi=0\), this vector field vanishes and the ordinary cosymplectic
Poisson picture is recovered.

The conformally deformed closedness conditions turn the bivector associated
with the odd pair into a Jacobi structure. More precisely, we have the following result, which can be proved in a manner similar to that used in \cite{AtEsLe23}:

\begin{proposition}\label{prop:LCCos-Jacobi-structure}
Let \((M,\tau,\Omega;\varphi)\) be an LCCos manifold satisfying
\[
d_\varphi\tau=0,
\qquad
d_\varphi\Omega=0,
\]
and let \(\Lambda^\tau\) be the bivector field associated with the
non-degenerate odd pair \((\tau,\Omega)\). Define the Lee vector field by
\[
\mathcal V_\varphi
=
\iota_\varphi\Lambda^\tau.
\]
Then the pair
\((\Lambda^\tau,\mathcal V_\varphi)\) determines a Jacobi structure on
\(M\). More precisely, it satisfies
\begin{equation}\label{ident-Jac+}
[\Lambda^\tau,\Lambda^\tau]_{SN}
=
-2\,\mathcal V_\varphi\wedge\Lambda^\tau,
\qquad
[\mathcal V_\varphi,\Lambda^\tau]_{SN}
=
0.
\end{equation}
\end{proposition}

Thus, the conformal deformation of the cosymplectic closedness conditions
transforms the Poisson bivector of ordinary cosymplectic geometry into the
Jacobi pair
\((\Lambda^\tau,\mathcal V_\varphi)\). In particular, when
the Lee form vanishes, one has \(\mathcal V_\varphi=0\), and
\eqref{ident-Jac+} reduces to the Poisson identity
\[
[\Lambda^\tau,\Lambda^\tau]_{SN}=0.
\]

Consequently, the LCCos bracket of two smooth functions
$F$ and $H$ in $C^\infty(M)$ is defined by
\begin{equation}\label{Pois-LCC}
\{F,H\}^{(\mathrm{LCCos})}
=
\Lambda^\tau(dF,dH)
+
F\,\mathcal{V}_{\varphi}(H)
-
H\,\mathcal{V}_{\varphi}(F).
\end{equation}
The identities in \eqref{ident-Jac+} imply that this bracket satisfies the Jacobi identity. In general, it is a first-order differential operator in each argument and does not satisfy the ordinary Leibniz rule. When \(\mathcal V_\varphi=0\), however, the additional terms vanish and the bracket reduces to a Poisson bracket.
 
\noindent
\textbf{From Local-to-Global.}
We now explain why the adjective ``locally conformally'' is appropriate for
the structure introduced above. Let \(M\) be a \((2n+1)\)-dimensional manifold
with an atlas \(\{U_\alpha\}\). Suppose that, on each chart \(U_\alpha\), one
is given a cosymplectic pair \((\tau_\alpha,\Omega_\alpha)\), namely
\begin{equation}\label{local-cosymplectic-conditions}
\tau_\alpha\wedge\Omega_\alpha^n\neq 0,\qquad d\tau_\alpha=0,
\quad
d\Omega_\alpha=0.
\end{equation}
The point is not to require these local pairs to coincide on intersections.
Rather, we ask that they become identical only after a local conformal change.
Thus, on every non-empty overlap \(U_\alpha\cap U_\beta\), we impose
\begin{equation}\label{lcc-loc-rel}
e^{\sigma_\alpha}\tau_\alpha
=
e^{\sigma_\beta}\tau_\beta,
\qquad
e^{\sigma_\alpha}\Omega_\alpha
=
e^{\sigma_\beta}\Omega_\beta,
\end{equation}
for some smooth functions $\sigma_\alpha$ in $C^\infty(U_\alpha)$.

This is the precise point at which our convention differs from the more
standard locally conformally cosymplectic convention.\footnote{In the
traditional formulation, the one-form and the two-form are usually assigned
different conformal weights: locally one writes
\(\tau|_{\alpha}=e^{\sigma_\alpha}\tau_\alpha\) and
\(\Omega|_{\alpha}=e^{2\sigma_\alpha}\Omega_\alpha\). Consequently, the
global equations become \(d\tau=\varphi\wedge\tau\) and
\(d\Omega=2\varphi\wedge\Omega\); see, for example,
\cite{AtEsLe23,ChLeMa91}.}
Here we deliberately use the same conformal weight for both local forms.
This choice is not merely cosmetic: it places \(\tau\) and \(\Omega\) under
the control of the same Lichnerowicz--de~Rham differential \(d_\varphi\).
Equivalently, both global structural equations will involve the same Lee
one-form with the same coefficient. This equal-weight convention is the one
that will be used below in the Hamiltonian and evolution dynamics.

Let us now pass from the local data to global forms. By
\eqref{lcc-loc-rel}, the locally defined forms
\begin{equation}\label{local-forms-cosymplectic}
\tau|_{\alpha}
=
e^{\sigma_\alpha}\tau_\alpha,
\qquad
\Omega|_{\alpha}
=
e^{\sigma_\alpha}\Omega_\alpha
\end{equation}
agree on overlaps \cite{AtEsLe23,ChLeMa91}. Hence they patch together to give a globally defined
one-form \(\tau\) and a globally defined two-form \(\Omega\) on \(M\).
Moreover, differentiating the overlap relations in \eqref{lcc-loc-rel} and
using \(d\tau_\alpha=0\) and \(d\Omega_\alpha=0\), one obtains
\(d\sigma_\alpha=d\sigma_\beta\) on \(U_\alpha\cap U_\beta\). Therefore, the
local exact one-forms \(d\sigma_\alpha\) determine a global closed one-form
\(\varphi\), namely
\begin{equation}\label{lcc-Lee-local}
\varphi|_{\alpha}
=
d\sigma_\alpha .
\end{equation}
This one-form is the Lee one-form of the locally conformal structure.

The non-degeneracy condition is also inherited from the local cosymplectic
data. Indeed, on each chart one has
\begin{equation}
\left.
\tau\wedge\Omega^n
\right|_{\alpha}
=
e^{(n+1)\sigma_\alpha}
\tau_\alpha\wedge\Omega_\alpha^n
\neq 0.
\end{equation}
Finally, the same local computation gives
\begin{equation}
\begin{split}
d\tau|_{\alpha}
&=
d\bigl(e^{\sigma_\alpha}\tau_\alpha\bigr)
=
d\sigma_\alpha\wedge e^{\sigma_\alpha}\tau_\alpha
=
\varphi|_{\alpha}\wedge\tau|_{\alpha},
\\
d\Omega|_{\alpha}
&=
d\bigl(e^{\sigma_\alpha}\Omega_\alpha\bigr)
=
d\sigma_\alpha\wedge e^{\sigma_\alpha}\Omega_\alpha
=
\varphi|_{\alpha}\wedge\Omega|_{\alpha}.
\end{split}
\end{equation} 
Thus, the global pair \((\tau,\Omega)\) satisfies
\begin{equation}\label{new-lcc-structure-equations}
d\tau
=
\varphi\wedge\tau,
\qquad
d\Omega
=
\varphi\wedge\Omega,
\qquad
\text{equivalently,}
\qquad
d_\varphi\tau=0,
\qquad
d_\varphi\Omega=0.
\end{equation}
The preceding construction yields the following local-to-global
characterization of the equal-weight LCCos structure.

\begin{proposition}
\label{prop:LCCos-local-to-global}
Let \(M\) be a \((2n+1)\)-dimensional manifold. The following statements are
equivalent:

\begin{enumerate}
\item The manifold \(M\) admits a non-degenerate odd pair
\((\tau,\Omega)\) and a closed one-form \(\varphi\) satisfying
\begin{equation}
\tau\wedge\Omega^n\neq0,
\qquad
d_\varphi\tau=0,
\qquad
d_\varphi\Omega=0.
\end{equation}

\item There exist an open cover \(\{U_\alpha\}\), local cosymplectic pairs
\((\tau_\alpha,\Omega_\alpha)\), and smooth functions
\(\sigma_\alpha\in C^\infty(U_\alpha)\) such that
\begin{equation}
d\tau_\alpha=0,
\qquad
d\Omega_\alpha=0,
\qquad
\tau_\alpha\wedge\Omega_\alpha^n\neq0,
\end{equation}
and, on every nonempty overlap \(U_\alpha\cap U_\beta\),
\begin{equation}
e^{\sigma_\alpha}\tau_\alpha
=
e^{\sigma_\beta}\tau_\beta,
\qquad
e^{\sigma_\alpha}\Omega_\alpha
=
e^{\sigma_\beta}\Omega_\beta,
\qquad
d\sigma_\alpha=d\sigma_\beta.
\end{equation}
\end{enumerate}

Under this correspondence, the global forms and the Lee form are locally
represented by
\begin{equation}
\tau|_\alpha
=
e^{\sigma_\alpha}\tau_\alpha,
\qquad
\Omega|_\alpha
=
e^{\sigma_\alpha}\Omega_\alpha,
\qquad
\varphi|_\alpha
=
d\sigma_\alpha.
\end{equation}
\end{proposition}

Hence, an equal-weight LCCos structure is locally conformal to an ordinary
cosymplectic structure, while its global geometry is encoded by the single
Lichnerowicz--de~Rham differential \(d_\varphi\).
 
\noindent
\textbf{A Generic LCCos Structure on \(\mathbb{R}\times T^*Q\).}
The extended cotangent bundle \(\mathbb{R}\times T^*Q\) provides a canonical
class of LCCos examples. Although one may start with a Lee form on the whole
space \(\mathbb{R}\times T^*Q\), we follow the natural choice discussed in
\cite{AtEsLe23} and take \(\varphi\) to be pulled back from the extended
configuration space \(\mathbb{R}\times Q\). Thus, by the same notation, let
\(\varphi\) denote a closed one-form on \(\mathbb{R}\times T^*Q\) with no
dependence on the fiber variables \(p_i\). We define
\begin{equation}\label{LCC-can}
\tau
=
d_{\varphi}s
=
ds-s\varphi,
\qquad
\Omega
=
-d_{\varphi}\Theta_Q
=
\Omega_Q+\varphi\wedge\Theta_Q,
\end{equation}
where \(\Theta_Q\) is the Liouville one-form and
\(\Omega_Q=-d\Theta_Q\) is the canonical symplectic form on \(T^*Q\).
Since \(d\varphi=0\), one immediately obtains
\begin{equation}
d\tau
=
\varphi\wedge\tau,
\qquad
d\Omega
=
\varphi\wedge\Omega.
\end{equation}
Hence
\((\mathbb{R}\times T^*Q,\tau,\Omega;\varphi)\) is an LCCos
manifold on the open region where the pair is non-degenerate.

On a chart \(U_\alpha\), let \((s,q^i,p_i)\) denote local coordinates. Then we have
\begin{equation}\label{varphi-coord}
\varphi|_{\alpha}
=
\phi\,ds+\varphi_i\,dq^i .
\end{equation}
Then \eqref{LCC-can} gives
\begin{equation}\label{LCC-RTQ-forms}
\tau|_{\alpha}
=
(1-s\phi)ds-s\varphi_i\,dq^i,
\qquad
\Omega|_{\alpha}
=
dq^i\wedge dp_i
-\frac{1}{2}C_{ij}\,dq^i\wedge dq^j
+
p_i\phi\,ds\wedge dq^i,
\end{equation}
where
\begin{equation}\label{Aij}
C_{ij}
=
\varphi_jp_i-\varphi_ip_j .
\end{equation}
On the non-degenerate region $1-s\phi\neq 0$, the corresponding Reeb vector field is given by
\begin{equation}\label{Reeb-LCC-RTQ}
\mathcal{R}^{\tau}|_{\alpha}
=
\frac{1}{1-s\phi}\frac{\partial}{\partial s}
+
\frac{p_i\phi}{1-s\phi}\frac{\partial}{\partial p_i}.
\end{equation} 

\subsection{Conformal Analysis of Cosymplectic Dynamics} 

In this subsection, we study dynamics on LCCos manifolds. The
purpose is to globalize the local cosymplectic Hamiltonian and evolution
equations by using the equal conformal weight convention. In this way, the
ordinary differential \(d\) is replaced by the Lichnerowicz--de~Rham
differential \(d_\varphi\), and the usual unit normalization in the evolution
equation is replaced by the conformal unit \(\mathbb I\).

\noindent \textbf{LCCos Hamiltonian Dynamics.} There are two natural approaches to the construction of LCCos Hamiltonian
dynamics. The first one directly deforms the ordinary cosymplectic
Hamiltonian equations. The second one begins with
ordinary cosymplectic Hamiltonian systems on the local charts and globalizes
them through the equal-weight conformal transition rules. We now present
these constructions separately.

\noindent
\textbf{First Approach: Direct Lichnerowicz--de~Rham Deformation.}
Let \((M,\tau,\Omega;\varphi)\) be an LCCos manifold and let
\(H\in C^\infty(M)\). Replacing \(dH\) in the ordinary cosymplectic
Hamiltonian equations by \(d_\varphi H\), we define the LCCos Hamiltonian
vector field \(X_H\) by
\begin{equation}\label{LCCos-Ham-global}
\iota_{X_H}\tau=0,
\qquad
\iota_{X_H}\Omega
=
d_\varphi H
-
\left\langle d_\varphi H,\mathcal R^\tau\right\rangle\tau.
\end{equation}
Equivalently, in terms of the bivector field \(\Lambda^\tau\), this vector
field is
\begin{equation}\label{LCCos-Ham-Jacobi-form} 
X_H
 =
\sharp_{\Lambda^\tau}(d_\varphi H)
 =
\sharp_{\Lambda^\tau}(dH-H\varphi)
=
\sharp_{\Lambda^\tau}(dH)
-
H\mathcal V_\varphi, 
\end{equation}
where $\mathcal V_\varphi
=
\iota_\varphi\Lambda^\tau.$ Thus, \(X_H\) is precisely the Hamiltonian vector field associated with the
Jacobi structure
\((\Lambda^\tau,\mathcal V_\varphi)\) established in
Proposition~\ref{prop:LCCos-Jacobi-structure}.

\noindent
\textbf{Second Approach: Globalization of Local Cosymplectic Dynamics.}
On each cosymplectic chart
\((U_\alpha,\tau_\alpha,\Omega_\alpha)\), let \(H_\alpha\) be a local
Hamiltonian function. The corresponding local cosymplectic Hamiltonian vector
field \(X_{H_\alpha}\) is determined by
\begin{equation}\label{localHam-LCCos}
\iota_{X_{H_\alpha}}\tau_\alpha=0,
\qquad
\iota_{X_{H_\alpha}}\Omega_\alpha
=
dH_\alpha
-
\mathcal R_\alpha(H_\alpha)\tau_\alpha.
\end{equation}
We assume that the local Hamiltonians transform with the same conformal
weight as the geometric data:
\begin{equation}\label{local-Hamiltonians-LCC}
e^{\sigma_\alpha}H_\alpha
=
e^{\sigma_\beta}H_\beta
\qquad
\text{on }
U_\alpha\cap U_\beta.
\end{equation}
They therefore define a global Hamiltonian function \(H\) through
\begin{equation}\label{global-H-LCC}
H|_\alpha
=
e^{\sigma_\alpha}H_\alpha.
\end{equation}
Since
\(\varphi|_\alpha=d\sigma_\alpha\), one obtains
\begin{equation}\label{LCCos-twisted-H-local}
\left.d_\varphi H\right|_\alpha
=
d\left(e^{\sigma_\alpha}H_\alpha\right)
-
d\sigma_\alpha\,e^{\sigma_\alpha}H_\alpha
=
e^{\sigma_\alpha}dH_\alpha.
\end{equation}

Under the equal-weight conformal transformations
\[
\tau|_\alpha=e^{\sigma_\alpha}\tau_\alpha,
\qquad
\Omega|_\alpha=e^{\sigma_\alpha}\Omega_\alpha,
\]
the Reeb vector fields satisfy
\begin{equation}\label{LCCos-Reeb-local-global}
\mathcal R^\tau|_\alpha
=
e^{-\sigma_\alpha}\mathcal R_\alpha.
\end{equation}
Consequently,
\begin{equation}
\left.
\left\langle d_\varphi H,\mathcal R^\tau\right\rangle
\right|_\alpha
=
\left\langle
e^{\sigma_\alpha}dH_\alpha,
e^{-\sigma_\alpha}\mathcal R_\alpha
\right\rangle
=
\mathcal R_\alpha(H_\alpha).
\end{equation}
Restricting the global equation \eqref{LCCos-Ham-global} to \(U_\alpha\)
therefore gives
\begin{equation}
\begin{split}
\left.
\left(
d_\varphi H
-
\left\langle d_\varphi H,\mathcal R^\tau\right\rangle\tau
\right)
\right|_\alpha
&=
e^{\sigma_\alpha}
\left(
dH_\alpha
-
\mathcal R_\alpha(H_\alpha)\tau_\alpha
\right)
\\
&=
e^{\sigma_\alpha}
\iota_{X_{H_\alpha}}\Omega_\alpha
\\
&=
\iota_{X_{H_\alpha}}\Omega|_\alpha.
\end{split}
\end{equation}
Similarly,
\[
\iota_{X_{H_\alpha}}\tau|_\alpha
=
e^{\sigma_\alpha}
\iota_{X_{H_\alpha}}\tau_\alpha
=
0.
\]
Thus, each \(X_{H_\alpha}\) satisfies the restriction of the global LCCos
Hamiltonian equations. By the non-degeneracy of the odd pair, the solution
of these equations is unique, and hence the local vector fields agree on
overlaps.

The equivalence of the two constructions can now be stated as follows.

\begin{proposition}
\label{prop:LCCos-Hamiltonian-equivalence}
Let \((M,\tau,\Omega;\varphi)\) be an LCCos manifold obtained from
equal-weight local cosymplectic data, and let the local Hamiltonians
\(\{H_\alpha\}\) satisfy \eqref{local-Hamiltonians-LCC}. Then the local
cosymplectic Hamiltonian vector fields \(X_{H_\alpha}\) glue directly to a
global vector field \(X_H\):
\begin{equation}\label{globalHam-LCCos}
X_H|_\alpha=X_{H_\alpha}.
\end{equation}
Moreover, this global vector field coincides with the vector field obtained
by the direct Lichnerowicz--de~Rham deformation:
\begin{equation}
\iota_{X_H}\tau=0,
\qquad
\iota_{X_H}\Omega
=
d_\varphi H
-
\left\langle d_\varphi H,\mathcal R^\tau\right\rangle\tau.
\end{equation}
Equivalently,
\begin{equation}
X_H
=
\sharp_{\Lambda^\tau}(dH)
-
H\mathcal V_\varphi.
\end{equation}
Thus, direct deformation and conformal globalization provide two equivalent
formulations of the same LCCos Hamiltonian dynamics.
\end{proposition}

The direct gluing in Proposition~\ref{prop:LCCos-Hamiltonian-equivalence} follows from the matching conformal weights of \(\Omega_\alpha\) and \(H_\alpha\), together with the compatibility of the local cosymplectic data. Assigning the same weight also to \(\tau_\alpha\) allows both structural equations to be expressed through the same Lichnerowicz--de~Rham differential \(d_\varphi\).

The Jacobi structure
\((\Lambda^\tau,\mathcal V_\varphi)\) also determines the algebraic and
geometric integrability of the resulting dynamics.

\begin{proposition}
\label{prop:LCCos-Hamiltonian-integrability}
Let \((M,\tau,\Omega;\varphi)\) be an LCCos manifold. For any
\(F,H\in C^\infty(M)\), the corresponding Hamiltonian vector fields satisfy
\begin{equation}\label{LCCos-Ham-Lie-closure}
[X_F,X_H]
=
-X_{\{F,H\}^{(\mathrm{LCCos})}},
\end{equation}
where
\begin{equation}
\{F,H\}^{(\mathrm{LCCos})}
=
\Lambda^\tau(dF,dH)
+
F\mathcal V_\varphi(H)
-
H\mathcal V_\varphi(F).
\end{equation}
\end{proposition}
Consequently, the LCCos Hamiltonian vector fields form an involutive family,
and the distribution
\begin{equation}
\mathcal D_{\mathrm{Ham}}
=
\operatorname{span}
\left\{
X_H:\ H\in C^\infty(M)
\right\}
\end{equation}
is integrable in the sense of Stefan--Sussmann.

\noindent \textbf{Coordinate Formulas for LCCos Hamiltonian Dynamics.}
In the local model \(\mathbb{R}\times T^*Q\), the LCCos Hamiltonian dynamics
can be written explicitly in the coordinates \((s,q^i,p_i)\). We use the Lee
form displayed in \eqref{varphi-coord} and the local formulas for
\(\tau\) and \(\Omega\) in \eqref{LCC-RTQ-forms}. For a Hamiltonian function
\(H=H(s,q,p)\), the LCCos Hamiltonian vector field \(X_H\) defined in
\eqref{LCCos-Ham-global} is computed to be
\begin{equation}\label{xiH-RTQ-local}
X_H
=
\frac{s\varphi_i}{1-s\phi}
\frac{\partial H}{\partial p_i}
\frac{\partial}{\partial s}
+
\frac{\partial H}{\partial p_i}
\frac{\partial}{\partial q^i}
+
\Bigl(
-\frac{\partial H}{\partial q^i}
+
\frac{C_{ij}}{1-s\phi}
\frac{\partial H}{\partial p_j}
+
\frac{\varphi_i}{1-s\phi}H
-
\frac{s\varphi_i}{1-s\phi}
\frac{\partial H}{\partial s}
\Bigr)
\frac{\partial}{\partial p_i}.
\end{equation}
Therefore, along an integral curve
\(t\mapsto (s(t),q^i(t),p_i(t))\), the LCCos
Hamiltonian equations of motion read
\begin{equation}\label{LCC-eom-RTQ}
\begin{split}
\frac{ds}{dt}
&=
\frac{s\varphi_i}{1-s\phi}
\frac{\partial H}{\partial p_i},
\\
\frac{dq^i}{dt}
&=
\frac{\partial H}{\partial p_i},
\\
\frac{dp_i}{dt}
&=
-\frac{\partial H}{\partial q^i}
+
\frac{C_{ij}}{1-s\phi}
\frac{\partial H}{\partial p_j}
+
\frac{\varphi_i}{1-s\phi}H
-
\frac{s\varphi_i}{1-s\phi}
\frac{\partial H}{\partial s}.
\end{split}
\end{equation}
These formulas display the effect of the Lee form directly. The \(q^i\)-part
keeps the usual canonical form, whereas the \(p_i\)-equation acquires
additional terms coming from the conformal deformation. More importantly, the
variable \(s\) is no longer frozen along the Hamiltonian flow. This happens
because the condition \(\iota_{X_H}\tau=0\) is imposed with
\(\tau=ds-s\varphi\), not with the ordinary form \(ds\). Hence, unlike in the
classical cosymplectic case, the Hamiltonian field may have a genuine
component in the \(s\)-direction.

In the particular case \(\phi=0\), the Lee form has no \(ds\)-component, that is 
$\varphi|_{\alpha}
=
\varphi_i\,dq^i$. 
Then \(1-s\phi=1\), and the LCCos Hamiltonian vector field
\eqref{xiH-RTQ-local} reduces to
\begin{equation}\label{xiH-RTQ-local-phi-zero}
X_H
=
s\varphi_i
\frac{\partial H}{\partial p_i}
\frac{\partial}{\partial s}
+
\frac{\partial H}{\partial p_i}
\frac{\partial}{\partial q^i}
+
\Bigl(
-\frac{\partial H}{\partial q^i}
+
C_{ij}
\frac{\partial H}{\partial p_j}
+
\varphi_i H
-
s\varphi_i
\frac{\partial H}{\partial s}
\Bigr)
\frac{\partial}{\partial p_i}.
\end{equation}
Accordingly, the equations of motion become
\begin{equation}\label{LCC-eom-RTQ-phi-zero}
\begin{split}
\frac{ds}{dt}
&=
s\varphi_i
\frac{\partial H}{\partial p_i},
\\
\frac{dq^i}{dt}
&=
\frac{\partial H}{\partial p_i},
\\
\frac{dp_i}{dt}
&=
-\frac{\partial H}{\partial q^i}
+
C_{ij}
\frac{\partial H}{\partial p_j}
+
\varphi_i H
-
s\varphi_i
\frac{\partial H}{\partial s}.
\end{split}
\end{equation}
Thus, even when \(\phi=0\), the dynamics is not necessarily the ordinary
cosymplectic dynamics. The remaining components \(\varphi_i\) of the Lee form
still produce conformal correction terms, and the variable \(s\) may still
evolve along the Hamiltonian flow. If, in addition, \(\varphi_i=0\), then the
Lee form vanishes, \(\varphi=0\), and one recovers the genuine cosymplectic
Hamiltonian dynamics. 

\noindent
\textbf{LCCos Evolution Dynamics.}
There are two natural approaches to LCCos evolution dynamics, parallel to
those used for the Hamiltonian dynamics. The first directly replaces the
ordinary exterior derivative by \(d_\varphi\) and the constant normalization
by the global conformal unit \(\mathbb I\). The second starts with local
cosymplectic evolution vector fields and globalizes them through the
equal-weight transition rules.

Throughout this discussion, we work in the exact-Lee sector in which the
global conformal unit \(\mathbb I\) exists. As established in
Subsection~\ref{sec-LdR}, this function is locally represented by
\begin{equation}\label{global-weighted-unit}
\mathbb I|_\alpha
=
e^{\sigma_\alpha}1_\alpha
\end{equation}
and satisfies
\begin{equation}
d_\varphi\mathbb I=0.
\end{equation}

\noindent
\textbf{First Approach: Direct Lichnerowicz--de~Rham Deformation.}
In ordinary cosymplectic evolution dynamics, the vector field is normalized
by the constant unit. In the conformally globalized setting, this
normalization is replaced by the \(d_\varphi\)-closed conformal unit
\(\mathbb I\). Accordingly, we define the LCCos evolution vector field
\(\mathcal E_H\) by
\begin{equation}\label{global-LCCos-evolution}
\iota_{\mathcal E_H}\tau
=
\mathbb I,
\qquad
\iota_{\mathcal E_H}\Omega
=
d_\varphi H
-
\left\langle d_\varphi H,\mathcal R^\tau\right\rangle\tau.
\end{equation}

The second equation in \eqref{global-LCCos-evolution} coincides with the
equation defining the LCCos Hamiltonian vector field \(X_H\). Since
\[
\iota_{\mathcal R^\tau}\tau=1,
\qquad
\iota_{\mathcal R^\tau}\Omega=0,
\]
the vector field
\begin{equation}\label{LCCos-evolution-Ham-Reeb}
\mathcal E_H
=
X_H+\mathbb I\,\mathcal R^\tau
\end{equation}
satisfies both equations in \eqref{global-LCCos-evolution}. Thus, the passage
from Hamiltonian dynamics to evolution dynamics is achieved by adding the
Reeb direction weighted by the global conformal unit.

\noindent
\textbf{Second Approach: Globalization of Local Evolution Dynamics.}
On each cosymplectic chart
\((U_\alpha,\tau_\alpha,\Omega_\alpha)\), let
\(\mathcal E_{H_\alpha}\) be the locally normalized cosymplectic evolution
vector field determined by
\begin{equation}\label{local-evolution-cos}
\iota_{\mathcal E_{H_\alpha}}\tau_\alpha
=
1_\alpha,
\qquad
\iota_{\mathcal E_{H_\alpha}}\Omega_\alpha
=
dH_\alpha
-
\mathcal R_\alpha(H_\alpha)\tau_\alpha,
\end{equation}
where \(1_\alpha\) is a nowhere-vanishing locally constant normalization.

The local normalizations are assumed to satisfy the weighted compatibility
condition
\begin{equation}\label{local-unit-gluing-evolution}
e^{\sigma_\alpha}1_\alpha
=
e^{\sigma_\beta}1_\beta
\qquad
\text{on }
U_\alpha\cap U_\beta.
\end{equation}
Similarly, the local Hamiltonians satisfy
\[
e^{\sigma_\alpha}H_\alpha
=
e^{\sigma_\beta}H_\beta.
\]
They therefore determine the global functions
\[
\mathbb I|_\alpha
=
e^{\sigma_\alpha}1_\alpha,
\qquad
H|_\alpha
=
e^{\sigma_\alpha}H_\alpha.
\]

We now rewrite the local equations in global variables. For the
normalization equation, one has
\begin{equation} 
\iota_{\mathcal E_{H_\alpha}}\tau|_\alpha
 =
\iota_{\mathcal E_{H_\alpha}}
\left(e^{\sigma_\alpha}\tau_\alpha\right)
 =
e^{\sigma_\alpha}
\iota_{\mathcal E_{H_\alpha}}\tau_\alpha
 =
e^{\sigma_\alpha}1_\alpha
 =
\mathbb I|_\alpha. 
\end{equation}
For the two-form equation, using
\[
\Omega|_\alpha=e^{\sigma_\alpha}\Omega_\alpha,
\qquad
\left.d_\varphi H\right|_\alpha
=
e^{\sigma_\alpha}dH_\alpha,
\qquad
\mathcal R^\tau|_\alpha
=
e^{-\sigma_\alpha}\mathcal R_\alpha,
\]
we obtain
\begin{equation}
\begin{split}
\iota_{\mathcal E_{H_\alpha}}\Omega|_\alpha
&=
e^{\sigma_\alpha}
\iota_{\mathcal E_{H_\alpha}}\Omega_\alpha
 =
e^{\sigma_\alpha}
\left(
dH_\alpha
-
\mathcal R_\alpha(H_\alpha)\tau_\alpha
\right)
\\
&=
\left.
\left(
d_\varphi H
-
\left\langle
d_\varphi H,\mathcal R^\tau
\right\rangle\tau
\right)
\right|_\alpha.
\end{split}
\end{equation}
Hence, every local evolution vector field
\(\mathcal E_{H_\alpha}\) satisfies the restriction of the global system
\eqref{global-LCCos-evolution}. By non-degeneracy, this system has a unique
solution, and the local evolution vector fields therefore agree on chart
overlaps.

The equivalence of the two approaches can be summarized as follows.

\begin{proposition}
\label{prop:LCCos-evolution-equivalence}
Let \((M,\tau,\Omega;\varphi)\) be an LCCos manifold in the exact-Lee sector,
and suppose that the local Hamiltonians and local constant normalizations
satisfy their equal-weight compatibility conditions. Then the local
cosymplectic evolution vector fields \(\mathcal E_{H_\alpha}\) glue directly
to a global vector field:
\begin{equation}\label{global-evolution-field}
\mathcal E_H|_\alpha
=
\mathcal E_{H_\alpha}.
\end{equation}
This global vector field coincides with the evolution vector field obtained
by direct Lichnerowicz--de~Rham deformation and is characterized by
\begin{equation}
\iota_{\mathcal E_H}\tau
=
\mathbb I,
\qquad
\iota_{\mathcal E_H}\Omega
=
d_\varphi H
-
\left\langle d_\varphi H,\mathcal R^\tau\right\rangle\tau.
\end{equation}
Equivalently, it satisfies
\begin{equation}
\mathcal E_H
=
X_H+\mathbb I\,\mathcal R^\tau.
\end{equation}
Thus, direct deformation and conformal globalization give two equivalent
formulations of LCCos evolution dynamics.
\end{proposition}

The conformal factor in the Reeb contribution is essential. Indeed, one may
also introduce the auxiliary vector field
\begin{equation}\label{auxiliary-LCCos-evolution}
\widetilde{\mathcal E}_H
=
X_H+\mathcal R^\tau,
\end{equation}
which satisfies
\begin{equation}
\iota_{\widetilde{\mathcal E}_H}\tau=1,
\qquad
\iota_{\widetilde{\mathcal E}_H}\Omega
=
d_\varphi H
-
\left\langle d_\varphi H,\mathcal R^\tau\right\rangle\tau.
\end{equation}
Although \(\widetilde{\mathcal E}_H\) is a globally defined vector field, it
does not arise from the equal-weight gluing of the locally normalized
evolution vector fields. The gluing-compatible evolution dynamics is
\(\mathcal E_H\), whose Reeb component is weighted by \(\mathbb I\). The two
vector fields coincide only when the conformal unit can be normalized to
\(\mathbb I=1\), which in turn gives \(\varphi=0\).

\noindent \textbf{Coordinate Formulas for LCCos Evolution Dynamics.}
The explicit form of the LCCos evolution field is obtained from the same local
data used in the Hamiltonian case. 
Fix a positive scalar conformal unit \(\mathbb I\). On each chart \(U_\alpha\), choose the Lee potential as \(\sigma_\alpha=\log(\mathbb I|_\alpha)\), thereby absorbing the local constant normalization into the potential. With this choice,
\begin{equation}
\varphi|_\alpha=d\sigma_\alpha,
\qquad
\mathbb I|_\alpha=e^{\sigma_\alpha(s,q)} .
\end{equation}
Hence, the conformal unit appears locally as the weighted function
\(e^{\sigma_\alpha(s,q)}\). Solving the defining equations of the LCCos
evolution vector field gives, for \(H=H(s,q,p)\),
\begin{equation}\label{EH-RTQ-local}
\begin{split}
\mathcal E_H \vert_\alpha
& =
\frac{1}{1-s\phi}
\left(
e^{\sigma_\alpha(s,q)}
+
s\varphi_i\frac{\partial H}{\partial p_i}
\right)
\frac{\partial}{\partial s}
+
\frac{\partial H}{\partial p_i}
\frac{\partial}{\partial q^i}
\\
&\qquad
+
\left(
-\frac{\partial H}{\partial q^i}
+
\frac{C_{ij}}{1-s\phi}
\frac{\partial H}{\partial p_j}
+
\frac{\varphi_i}{1-s\phi}H
-
\frac{s\varphi_i}{1-s\phi}
\frac{\partial H}{\partial s}
+
\frac{p_i\phi}{1-s\phi}
e^{\sigma_\alpha(s,q)}
\right)
\frac{\partial}{\partial p_i}.
\end{split}
\end{equation}
Therefore, one obtains
\begin{equation}\label{LCCos-evo-eq}
\begin{split}
\frac{ds}{dt}
&=
\frac{1}{1-s\phi}
\left(
e^{\sigma_\alpha(s,q)}
+
s\varphi_i\frac{\partial H}{\partial p_i}
\right),
\\
\frac{dq^i}{dt}
&=
\frac{\partial H}{\partial p_i},
\\
\frac{dp_i}{dt}
&=
-\frac{\partial H}{\partial q^i}
+
\frac{C_{ij}}{1-s\phi}
\frac{\partial H}{\partial p_j}
+
\frac{\varphi_i}{1-s\phi}H
-
\frac{s\varphi_i}{1-s\phi}
\frac{\partial H}{\partial s}
+
\frac{p_i\phi}{1-s\phi}
e^{\sigma_\alpha(s,q)} .
\end{split}
\end{equation}
This system differs from the Hamiltonian system in \eqref{LCC-eom-RTQ} by the
presence of the conformal unit in the Reeb direction. The reason is that the
evolution field is normalized by
\(\iota_{\mathcal E_H}\tau=\mathbb I\), whereas the Hamiltonian field satisfies
\(\iota_{X_H}\tau=0\). Consequently, the equation for \(ds/dt\)
contains the additional term \(e^{\sigma_\alpha(s,q)}\). Likewise, the last
term in the \(p_i\)-equation is precisely the local contribution of the
weighted Reeb component \(\mathbb I\mathcal R^\tau\).

In the particular case \(\phi=0\), the Lee form has no \(ds\)-component,
namely 
$\varphi|_\alpha
=
\varphi_i\,dq^i$. 
Equivalently, the local conformal factor \(\sigma_\alpha\) is independent of
\(s\). Hence \(1-s\phi=1\), and the evolution vector field
\eqref{EH-RTQ-local} reduces to
\begin{equation}\label{EH-RTQ-local-phi-zero} 
\mathcal E_H \vert_\alpha
 =
\left(
e^{\sigma_\alpha(q)}
+
s\varphi_i\frac{\partial H}{\partial p_i}
\right)
\frac{\partial}{\partial s}
+
\frac{\partial H}{\partial p_i}
\frac{\partial}{\partial q^i}
+
\left(
-\frac{\partial H}{\partial q^i}
+
C_{ij}
\frac{\partial H}{\partial p_j}
+
\varphi_i H
-
s\varphi_i
\frac{\partial H}{\partial s}
\right)
\frac{\partial}{\partial p_i}. 
\end{equation}
Accordingly, the evolution equations become
\begin{equation}\label{LCCos-evo-eq-phi-zero}
\begin{split}
\frac{ds}{dt}
&=
e^{\sigma_\alpha(q)}
+
s\varphi_i\frac{\partial H}{\partial p_i},
\\
\frac{dq^i}{dt}
&=
\frac{\partial H}{\partial p_i},
\\
\frac{dp_i}{dt}
&=
-\frac{\partial H}{\partial q^i}
+
C_{ij}
\frac{\partial H}{\partial p_j}
+
\varphi_i H
-
s\varphi_i
\frac{\partial H}{\partial s}.
\end{split}
\end{equation}
Thus, even when \(\phi=0\), the system is still genuinely locally conformal
unless the remaining components \(\varphi_i\) also vanish. If, in addition,
\(\varphi_i=0\), then the Lee form vanishes, \(\varphi=0\), the conformal unit
may be normalized to \(\mathbb I=1\), and one recovers the genuine
cosymplectic evolution dynamics.
 
\subsection{Time-Dependent Mass Dynamics in an LCCos Environment}
\label{sec:LCCos-mass-example}

In this subsection, we present a mechanical model illustrating how the LCCos
framework modifies time-dependent Hamiltonian dynamics. The example is based
on a particle whose kinetic energy is determined by a possibly
time-dependent Riemannian metric, or equivalently by a time-dependent mass
tensor. From a mechanical viewpoint, such a metric describes a system whose
kinetic-energy geometry, and hence its effective mass or inertia distribution,
changes explicitly with time. In the absence of the Lee form, the model
reduces to the ordinary cosymplectic description of non-autonomous
Hamiltonian mechanics. When the Lee form is nonzero, additional conformal
force terms appear in the equations of motion.

Let $Q$ be a three-dimensional manifold with local coordinates
$\mathbf q=(q^1,q^2,q^3)$. The corresponding cosymplectic phase space is the
seven-dimensional manifold $\mathbb R\times T^*Q$, endowed with induced local
coordinates $(s,\mathbf q,\mathbf p)$. Consider the time-dependent
Hamiltonian function
\begin{equation}\label{LCCos-time-dependent-Hamiltonian}
H(s,\mathbf q,\mathbf p)
=
\frac{1}{2}\mathbf p\cdot G(s,\mathbf q)\mathbf p
+
V(s,\mathbf q).
\end{equation}
Here $G(s,\mathbf q)=[G^{ij}(s,\mathbf q)]$ denotes the inverse matrix of a
time-dependent Riemannian metric
$\mathcal G(s,\mathbf q)=[\mathcal G_{ij}(s,\mathbf q)]$ on $Q$, and
$V(s,\mathbf q)$ is a time-dependent potential. Equivalently,
$\mathcal G(s,\mathbf q)$ may be interpreted as a time-dependent mass tensor,
while $G(s,\mathbf q)$ represents the corresponding inverse mass tensor.

To isolate the locally conformal deformation from the explicit
$s$-dependence of the mechanical data, we assume that the Lee form has no
$ds$-component. Thus, we take
\begin{equation}\label{LCCos-Lee-q-only}
\varphi
=
\boldsymbol{\varphi}(\mathbf q)\cdot d\mathbf q,
\end{equation}
where
$\boldsymbol{\varphi}=(\varphi_1,\varphi_2,\varphi_3)$. Since the Lee form is
closed, its components satisfy
\begin{equation}
\partial_i\varphi_j-\partial_j\varphi_i=0.
\end{equation}
Note that the assumption \eqref{LCCos-Lee-q-only} does not remove the
explicit $s$-dependence of the mechanical system: both
$\mathcal G(s,\mathbf q)$ and $V(s,\mathbf q)$ are still allowed to depend
explicitly on $s$.

For this choice of Lee form, the general LCCos structure simplifies
considerably. In particular, the Reeb vector field becomes
\begin{equation}
\mathcal R^\tau
=
\frac{\partial}{\partial s}.
\end{equation}
Moreover, the skew-symmetric tensor entering the local LCCos equations is
\begin{equation}
C_{ij}
=
\varphi_jp_i-\varphi_ip_j.
\end{equation}
In three dimensions, we use the cross-product notation as a convenient
shorthand for this component expression. More precisely, for a vector
$\mathbf a$,
\begin{equation}\label{LCCos-C-cross-product}
(C\mathbf a)_i
=
C_{ij}a^j,
\qquad
C\mathbf a
=
(\boldsymbol{\varphi}\times\mathbf p)\times\mathbf a.
\end{equation}

\noindent
\textbf{LCCos Hamiltonian Dynamics.}
Substituting the Hamiltonian
\eqref{LCCos-time-dependent-Hamiltonian} and the Lee form
\eqref{LCCos-Lee-q-only} into the general LCCos Hamiltonian equations
\eqref{LCC-eom-RTQ}, we obtain
\begin{equation}\label{LCCos-mass-Hamilton-equations}
\begin{split}
\frac{ds}{dt}
&=
s\,\boldsymbol{\varphi}\cdot G\mathbf p,
\\
\frac{d\mathbf q}{dt}
&=
G\mathbf p,
\\
\frac{d\mathbf p}{dt}
&=
-
\nabla_{\mathbf q}H
+
(\boldsymbol{\varphi}\times\mathbf p)\times G\mathbf p
+
\boldsymbol{\varphi}H
-
s\boldsymbol{\varphi}
\frac{\partial H}{\partial s}.
\end{split}
\end{equation}
Here
\begin{equation}\label{LCCos-H-derivatives}
\nabla_{\mathbf q}H
=
\frac{1}{2}
\nabla_{\mathbf q}
\bigl(
\mathbf p\cdot G(s,\mathbf q)\mathbf p
\bigr)
+
\nabla_{\mathbf q}V(s,\mathbf q),
\end{equation}
whereas
\begin{equation}
\frac{\partial H}{\partial s}
=
\frac{1}{2}
\mathbf p\cdot
\frac{\partial G}{\partial s}(s,\mathbf q)\mathbf p
+
\frac{\partial V}{\partial s}(s,\mathbf q).
\end{equation}

Thus, the configuration equation retains the canonical form
$\dot{\mathbf q}=G\mathbf p$, whereas the momentum equation acquires
additional terms determined by the Lee form. In particular,
$(\boldsymbol{\varphi}\times\mathbf p)\times G\mathbf p$ is the
three-dimensional representation of the skew-symmetric contribution
$C_{ij}G^{jk}p_k$. The remaining terms in the Lee direction involve both the
Hamiltonian itself and its explicit $s$-dependence.

An additional feature of the LCCos Hamiltonian dynamics is that the variable
$s$ is no longer necessarily constant along the Hamiltonian flow. Indeed,
using $\dot{\mathbf q}=G\mathbf p$, the first equation in
\eqref{LCCos-mass-Hamilton-equations} becomes
\begin{equation}\label{LCCos-Hamiltonian-s-dot}
\dot s
=
s\,\boldsymbol{\varphi}\cdot\dot{\mathbf q}.
\end{equation}
Hence, although the Lee form has no $ds$-component, its configurational
component couples the Hamiltonian motion to the $s$-direction.

We now pass to the corresponding Newton form. Since
\begin{equation}
\dot{\mathbf q}
=
G\mathbf p,
\qquad
\mathbf p
=
\mathcal G\dot{\mathbf q},
\end{equation}
the $s$-derivative of the Hamiltonian, when expressed in velocity variables,
takes the form
\begin{equation}\label{LCCos-Hs-velocity}
\frac{\partial H}{\partial s}
=
-
\frac{1}{2}
\dot{\mathbf q}\cdot
(\partial_s\mathcal G)
\dot{\mathbf q}
+
\frac{\partial V}{\partial s}.
\end{equation}
Accordingly, the LCCos Hamiltonian dynamics can be written as
\begin{equation}\label{LCCos-Hamiltonian-Newton-form}
\ddot{\mathbf q}
+
\boldsymbol{\Gamma}_{s}
(\dot{\mathbf q},\dot{\mathbf q})
=
-
G\nabla_{\mathbf q}V
-
G\bigl(
(\partial_s\mathcal G)\dot{\mathbf q}
\bigr)\dot s
+
G\mathbf u
+
\left(
\frac{s}{2}
\dot{\mathbf q}\cdot
(\partial_s\mathcal G)
\dot{\mathbf q}
-
s\frac{\partial V}{\partial s}
\right)
G\boldsymbol{\varphi},
\end{equation}
together with
\begin{equation}
\dot s
=
s\,\boldsymbol{\varphi}\cdot\dot{\mathbf q}.
\end{equation}
Here,
$\boldsymbol{\Gamma}_{s}(\dot{\mathbf q},\dot{\mathbf q})$ denotes the
Christoffel contribution associated with the metric
$\mathcal G(s,\mathbf q)$ at fixed $s$, whose $\ell$-th component is
\begin{equation}
\bigl(
\boldsymbol{\Gamma}_{s}
(\dot{\mathbf q},\dot{\mathbf q})
\bigr)^\ell
=
\Gamma^\ell_{ij}(s,\mathbf q)
\dot q^i\dot q^j.
\end{equation}
The locally conformal contribution $\mathbf u$ is defined by
\begin{equation}\label{LCCos-u-term}
\mathbf u
=
(\boldsymbol{\varphi}\times
\mathcal G\dot{\mathbf q})
\times\dot{\mathbf q}
+
\boldsymbol{\varphi}
\left(
\frac{1}{2}
\dot{\mathbf q}\cdot
\mathcal G\dot{\mathbf q}
+
V(s,\mathbf q)
\right).
\end{equation}

The vector $\mathbf u$ collects the Lee-form contributions that are already
present independently of the explicit $s$-derivative of the Hamiltonian.
Its first term is the three-dimensional skew contribution produced by the
conformal twisting, whereas the second term is an energy-weighted
contribution in the Lee direction. Since
$\boldsymbol{\varphi}$ originates from a one-form, $\mathbf u$ is naturally
interpreted in covector form, while the inverse metric $G$ raises the index.
Consequently, $G\mathbf u$ represents the corresponding acceleration and
measures the deviation of the LCCos dynamics from the ordinary cosymplectic
mechanical motion.

The remaining Lee-direction term in
\eqref{LCCos-Hamiltonian-Newton-form},
\begin{equation}\label{LCCos-explicit-time-Lee-force}
\left(
\frac{s}{2}
\dot{\mathbf q}\cdot
(\partial_s\mathcal G)
\dot{\mathbf q}
-
s\frac{\partial V}{\partial s}
\right)
G\boldsymbol{\varphi},
\end{equation}
is generated specifically by the explicit $s$-dependence of the
Hamiltonian. Indeed, it originates from
$-s\boldsymbol{\varphi}\,\partial_sH$ in
\eqref{LCCos-mass-Hamilton-equations} after using
\eqref{LCCos-Hs-velocity}. Thus, this contribution describes how the explicit
time dependence of the mass tensor and the potential is transmitted to the
mechanical dynamics through the Lee direction.

Notice also that the term
\begin{equation}
-
G\bigl(
(\partial_s\mathcal G)\dot{\mathbf q}
\bigr)\dot s
\end{equation}
is present because the metric itself varies with $s$. In the LCCos
Hamiltonian dynamics, this effect is coupled to the conformal geometry
through
$\dot s=s\,\boldsymbol{\varphi}\cdot\dot{\mathbf q}$.

\noindent
\textbf{LCCos Evolution Dynamics.}
We now turn to the corresponding evolution dynamics. As in the general
construction, we work in the exact-Lee sector in which the global conformal
unit exists. Since the Lee form has only configurational components, we may
locally write
\begin{equation}
\varphi=d\sigma,
\qquad
\sigma=\sigma(\mathbf q),
\end{equation}
and hence
\begin{equation}
\mathbb I
=
e^{\sigma(\mathbf q)}.
\end{equation}

The LCCos evolution vector field is obtained from the Hamiltonian vector
field by adding the weighted Reeb contribution,
\begin{equation}
\mathcal E_H
=
X_H
+
\mathbb I\,\mathcal R^\tau.
\end{equation}
Since $\mathcal R^\tau=\partial/\partial s$, the Reeb contribution affects
only the $s$-component of the evolution vector field. Therefore, the local
evolution equations reduce to
\begin{equation}\label{LCCos-evolution-coordinate-equations}
\begin{split}
\dot s_{\mathrm{evo}}
&=
e^{\sigma(\mathbf q)}
+
s\,\boldsymbol{\varphi}\cdot\dot{\mathbf q},
\\
\dot{\mathbf q}
&=
G\mathbf p,
\\
\dot{\mathbf p}
&=
\left(
\dot{\mathbf p}
\right)_{\mathrm{Ham}},
\end{split}
\end{equation}
where
$\left(\dot{\mathbf p}\right)_{\mathrm{Ham}}$ denotes the momentum equation
in \eqref{LCCos-mass-Hamilton-equations}.

Thus, under the assumption that the Lee form has no $ds$-component, the
weighted Reeb field produces no direct correction in the momentum equation.
Its effect enters the Newton dynamics through the evolution of $s$, since the
metric and the potential depend explicitly on this variable.

Accordingly, the Newton form of the LCCos evolution dynamics is
\begin{equation}\label{LCCos-evolution-Newton-form}
\begin{split}
\ddot{\mathbf q}
+
\boldsymbol{\Gamma}_{s}
(\dot{\mathbf q},\dot{\mathbf q})
={}&
-
G\nabla_{\mathbf q}V
-
G\bigl(
(\partial_s\mathcal G)\dot{\mathbf q}
\bigr)
\dot s_{\mathrm{evo}}
+
G\mathbf u
\\
&+
\left(
\frac{s}{2}
\dot{\mathbf q}\cdot
(\partial_s\mathcal G)
\dot{\mathbf q}
-
s\frac{\partial V}{\partial s}
\right)
G\boldsymbol{\varphi},
\end{split}
\end{equation}
where
\begin{equation}\label{LCCos-evolution-s-dot}
\dot s_{\mathrm{evo}}
=
e^{\sigma(\mathbf q)}
+
s\,\boldsymbol{\varphi}\cdot\dot{\mathbf q}.
\end{equation}

Comparing \eqref{LCCos-evolution-Newton-form} with the Hamiltonian Newton
equation \eqref{LCCos-Hamiltonian-Newton-form}, the difference between the
two acceleration equations is therefore entirely due to the weighted Reeb
contribution to the evolution of $s$. More precisely,
\begin{equation}\label{LCCos-evolution-extra-force}
-
e^{\sigma(\mathbf q)}
G\bigl(
(\partial_s\mathcal G)\dot{\mathbf q}
\bigr)
\end{equation}
is the additional term obtained when passing from the Hamiltonian dynamics
to the evolution dynamics. This term vanishes when the mass metric is
independent of $s$.

Finally, if the Lee form vanishes,
$\boldsymbol{\varphi}=0$, then
\begin{equation}
\dot s_{\mathrm{Ham}}=0,
\end{equation}
and the LCCos Hamiltonian dynamics reduces to the ordinary cosymplectic
Hamiltonian dynamics on fixed-$s$ slices. In the evolution picture, after
choosing the standard normalization $\mathbb I=1$, one obtains
\begin{equation}
\dot s_{\mathrm{evo}}=1,
\end{equation}
and hence the usual non-autonomous cosymplectic evolution dynamics is
recovered.

\section{Locally Conformal Analysis of Contact Framework}\label{LC-Sec-12} 

 In this section, we begin with coorientable contact manifolds in light of the odd pair discussion in Subsection~\ref{sec:ndop}. 
We then carry the locally conformal analysis to the contact
framework. Unlike the cosymplectic case, contact geometry already possesses an
intrinsic conformal freedom, since multiplying a contact form by a nowhere
vanishing function preserves the underlying contact distribution. Nevertheless,
this conformal flexibility does not automatically describe how local contact
Hamiltonian systems should be glued into a global dynamics. Our aim is to make
this point precise by introducing a locally conformal contact structure through
the Lichnerowicz--de~Rham differential. In this way, the ordinary relation
\(\Omega=d\eta\) is replaced by the twisted relation
\(\Omega=d_{\varphi}\eta\), and the Lee form records the conformal mismatch
between local contact descriptions. We then study the resulting Jacobi
structure, derive the associated Hamiltonian equations, and show how the
locally conformal terms modify the usual contact dissipative dynamics.

\subsection{Jacobi Character for an Odd Pair: Contact Geometry}\label{sec:Jacko}

We once again start with a non-degenerate odd manifold \((M,\Omega,\eta)\). In order to
be consistent with the standard contact terminology, we use the notation
\(\eta\) for the one-form and write the defining pair in the order
\((\Omega,\eta)\). The non-degeneracy condition is 
$\Omega^n\wedge\eta\neq 0$. 
As before, this pair determines the Reeb vector field \(\mathcal R^\eta\) and
the bivector field \(\Lambda^\eta\). More precisely,
\begin{equation}\label{Reeb-eta}
\iota_{\mathcal R^\eta}\eta=1,
\qquad
\iota_{\mathcal R^\eta}\Omega=0,
\end{equation}
and
\begin{equation}\label{Lambda-eta}
\Lambda^\eta(\alpha,\beta)
=
\Omega(\sharp^\eta\alpha,\sharp^\eta\beta).
\end{equation}
Thus, by the equivalence between the covariant and contravariant descriptions,
we may also represent the same geometric data by the pair
\((\Lambda^\eta,\mathcal R^\eta)\).

The essential difference in the present subsection is that the bracket on
functions is no longer determined only by the bivector field. Instead, it is
defined by using both the bivector field \(\Lambda^\eta\) and the Reeb vector
field \(\mathcal R^\eta\). Referring to the contravariant pair
\((\Lambda^\eta,\mathcal R^\eta)\), we define a skew-symmetric bracket on
\(C^\infty(M)\) by
\begin{equation}\label{almost-jacobi-bracket-eta}
\{F,H\}^{\mathcal R}
=
\Lambda^\eta(dF,dH)
+
F\mathcal R^\eta(H)
-
H\mathcal R^\eta(F).
\end{equation}
Compared with the almost Poisson bracket
\eqref{almost-poisson-bracket-tau}, this expression contains two additional
terms determined by the Reeb vector field. These terms are responsible for
the characteristic first-order nature of the Jacobi-type bracket.

In general, this bracket satisfies neither the Leibniz identity nor the Jacobi
identity. Instead, it satisfies the first-order differential property
\begin{equation}\label{almost-jacobi-first-order-eta}
\{F,KH\}^{\mathcal R}
=
K\{F,H\}^{\mathcal R}
+
H\{F,K\}^{\mathcal R}
+
KH\mathcal R^\eta(F).
\end{equation}
Hence, we refer to \eqref{almost-jacobi-bracket-eta} as an almost Jacobi
bracket. At this stage, the Jacobi identity is not guaranteed.

The Hamiltonian vector field \cite{Bravetti2017,de2019contact,leonsardon1,EsLaLeMa21,Simoes2020b} associated with the almost Jacobi bracket is
defined by
\begin{equation}\label{almost-jacobi-Ham-vf-eta}
X_H^{\mathcal R}
=
\sharp_{\Lambda^\eta}(dH)
-
H\mathcal R^\eta
=
\sharp^\eta(dH)
-
\bigl(\mathcal R^\eta(H)+H\bigr)\mathcal R^\eta .
\end{equation}
Equivalently, it is characterized intrinsically by
\begin{equation}\label{almost-jacobi-Ham-vf-char-eta}
\iota_{X_H^{\mathcal R}}\Omega
=
dH-\mathcal R^\eta(H)\eta,
\qquad
\iota_{X_H^{\mathcal R}}\eta
=
-H.
\end{equation}
The second identity shows the essential difference from the almost Poisson
case: the Hamiltonian vector field is no longer tangent to the distribution
\(\ker\eta\), except on the zero level set of \(H\).

In contrast with the almost Poisson case, the Hamiltonian function is not
conserved in general. Indeed,
\begin{equation}\label{Hamiltonian-evolution-XHR-eta}
X_H^{\mathcal R}(H)
=
-H\mathcal R^\eta(H).
\end{equation}
Consequently, the zero level set \(H^{-1}(0)\) is invariant under the flow of
\(X_H^{\mathcal R}\).

Let
\begin{equation}
\mathrm{Vol}=\Omega^n\wedge\eta
\end{equation}
be the volume form induced by the non-degenerate odd pair. Using Cartan's
formula together with the defining relations
\eqref{almost-jacobi-Ham-vf-char-eta}, one obtains
\begin{equation}\label{volume-XHR-eta}
\begin{split}
\mathcal L_{X_H^{\mathcal R}}\mathrm{Vol}
={}&
(\iota_{X_H^{\mathcal R}}d\eta)\wedge\Omega^n
+
n\eta\wedge(\iota_{X_H^{\mathcal R}}d\Omega)\wedge\Omega^{n-1}
-
\mathcal R^\eta(H)\mathrm{Vol}
-
n\mathcal R^\eta(H)\eta\wedge d\eta\wedge\Omega^{n-1}.
\end{split}
\end{equation}
Thus, the volume form is not preserved in general under the flow of
\(X_H^{\mathcal R}\).

Moreover, the almost Jacobi Hamiltonian vector field acts on functions as
\begin{equation}\label{almost-jacobi-Ham-action-eta}
X_H^{\mathcal R}(F)
=
\{F,H\}^{\mathcal R}
-
F\mathcal R^\eta(H).
\end{equation}
In particular, on the zero level set \(H^{-1}(0)\), the Reeb component
\(-H\mathcal R^\eta\) vanishes. Hence,
\begin{equation}\label{zero-level-XHR-XH-eta}
X_H^{\mathcal R}\big|_{H=0}
=
\sharp_{\Lambda^\eta}(dH)\big|_{H=0}.
\end{equation}
 
The possible failure of almost Jacobi Hamiltonian vector fields to close
under the Jacobi--Lie bracket is governed by the two tensors appearing in
the Jacobi conditions. The preceding computation can be summarized as
follows.

\begin{proposition}
\label{prop:almost-Jacobi-Hamiltonian-obstruction}
Let \(M\) be equipped with a non-degenerate odd pair
\((\Omega,\eta)\), and let
\((\Lambda^\eta,\mathcal R^\eta)\) be its associated contravariant pair.
Then, for every \(F,H,K\in C^\infty(M)\), the corresponding almost Jacobi
Hamiltonian vector fields satisfy
\begin{equation}\label{Jac-integ-almost-jacobi-eta}
\begin{split}
[X_F^{\mathcal R},X_H^{\mathcal R}](K)
={}&
-
X_{\{F,H\}^{\mathcal R}}^{\mathcal R}(K)
+
\frac{1}{2}
\left(
[\Lambda^\eta,\Lambda^\eta]_{SN}
+
2\mathcal R^\eta\wedge\Lambda^\eta
\right)
(dF,dH,dK)
\\
&+
H\bigl(\mathcal L_{\mathcal R^\eta}\Lambda^\eta\bigr)(dK,dF)
-
F\bigl(\mathcal L_{\mathcal R^\eta}\Lambda^\eta\bigr)(dK,dH).
\end{split}
\end{equation}
Consequently, the Hamiltonian assignment
\[
H\longmapsto X_H^{\mathcal R}
\]
satisfies
\begin{equation}
[X_F^{\mathcal R},X_H^{\mathcal R}]
=
-
X_{\{F,H\}^{\mathcal R}}^{\mathcal R}
\end{equation}
for all \(F,H\in C^\infty(M)\) if and only if
\begin{equation}\label{almost-Jacobi-obstruction-vanishing}
[\Lambda^\eta,\Lambda^\eta]_{SN}
=
-2\mathcal R^\eta\wedge\Lambda^\eta,
\qquad
[\mathcal R^\eta,\Lambda^\eta]_{SN}
=
\mathcal L_{\mathcal R^\eta}\Lambda^\eta
=
0.
\end{equation}
\end{proposition}

Thus, the two terms 
$[\Lambda^\eta,\Lambda^\eta]_{SN}
+
2\mathcal R^\eta\wedge\Lambda^\eta$
 and 
$\mathcal L_{\mathcal R^\eta}\Lambda^\eta$ 
measure precisely the failure of the almost Jacobi pair to satisfy the
Jacobi conditions. As shown below, both obstructions vanish in the contact
case.

\noindent
\textbf{Coorientable Contact Manifolds.}
A distinguished case is obtained when the two-form \(\Omega\) is exact and is
given by the exterior derivative of \(\eta\), namely \(\Omega=d\eta\). Then
the non-degeneracy condition becomes
\begin{equation}\label{cont-man-def}
\Omega^n\wedge\eta\neq 0,
\qquad
\Omega=d\eta.
\end{equation}
In this case, \(\eta\) is called a contact form, and
\((M,\Omega,\eta)\), with \(\Omega=d\eta\), is called a coorientable contact
manifold. Since the two-form \(\Omega\) is completely determined by \(\eta\),
one often denotes the same contact manifold simply by \((M,\eta)\).

Geometrically, the contact one-form \(\eta\) determines the hyperplane
distribution
\begin{equation}\label{contact-distribution-eta}
\mathcal D
=
\ker\eta
\subset TM.
\end{equation}
The non-degeneracy condition \((d\eta)^n\wedge\eta\neq 0\) is equivalent to
the maximal non-integrability of this distribution. The distribution
\(\mathcal D\) is called the contact structure. Note that, for a fixed
contact structure, the contact form is not unique. Indeed, if \(\eta\) is a
contact form and \(u\) is a nowhere-vanishing smooth function on \(M\), then
\(u\eta\) defines the same contact distribution, since
\begin{equation}\label{same-contact-distribution}
\ker(u\eta)=\ker\eta.
\end{equation}
Moreover, \(u\eta\) is again a contact form, because
\begin{equation}\label{rescaled-contact-form}
(u\eta)\wedge d(u\eta)^n
=
u^{n+1}(d\eta)^n\wedge\eta
\neq 0.
\end{equation}
 
The obstructions identified in
Proposition~\ref{prop:almost-Jacobi-Hamiltonian-obstruction} vanish when the
non-degenerate odd pair is induced by a contact form. This gives the
following result.

\begin{proposition}\label{prop:contact-Jacobi-structure}
Let \(M\) be a \((2n+1)\)-dimensional contact manifold with contact form
\(\eta\), associated two-form \(\Omega=d\eta\), Reeb vector field
\(\mathcal R^\eta\), and bivector field \(\Lambda^\eta\). Then the pair
\((\Lambda^\eta,\mathcal R^\eta)\) satisfies the Jacobi identities
\begin{equation}\label{contact-Jacobi-conditions}
[\Lambda^\eta,\Lambda^\eta]_{SN}
=
-2\mathcal R^\eta\wedge\Lambda^\eta,
\qquad
[\mathcal R^\eta,\Lambda^\eta]_{SN}
=
0.
\end{equation}
Consequently, the bracket
\begin{equation}
\{F,H\}^{\mathcal R}
=
\Lambda^\eta(dF,dH)
+
F\mathcal R^\eta(H)
-
H\mathcal R^\eta(F)
\end{equation}
satisfies the Jacobi identity and therefore defines a genuine Jacobi bracket
on \(C^\infty(M)\).
\end{proposition}

It follows from
Proposition~\ref{prop:almost-Jacobi-Hamiltonian-obstruction} that the defect
terms in \eqref{Jac-integ-almost-jacobi-eta} vanish identically. Therefore,
the contact Hamiltonian vector fields satisfy
\begin{equation}\label{Jac-integ}
[X_F^{\mathcal R},X_H^{\mathcal R}]
=
-
X_{\{F,H\}^{\mathcal R}}^{\mathcal R}.
\end{equation}
Hence, the space
\begin{equation}
\mathfrak X_{\mathrm{Ham}}^{\mathcal R}(M)
=
\left\{
X_H^{\mathcal R}:\ H\in C^\infty(M)
\right\}
\end{equation}
is closed under the Jacobi--Lie bracket and therefore forms a Lie algebra.
Thus, the contact closedness condition \(\Omega=d\eta\) removes the
obstruction present for a general non-degenerate odd pair and turns the
almost Jacobi Hamiltonian dynamics into a genuine Jacobi Hamiltonian
dynamics. 

This closure has the usual integrability consequence. Namely, the
characteristic distribution generated by Jacobi Hamiltonian vector fields is
involutive and hence integrable by the Frobenius theorem. In the contact case,
this distribution is in fact the whole tangent bundle, since
\(\operatorname{Im}\sharp_{\Lambda^\eta}=\ker\eta\) and
\(\mathcal R^\eta\) is transverse to \(\ker\eta\). Consequently, the leaves of the characteristic foliation are the connected components of \(M\). In this sense, Jacobi geometry preserves the integrability aspect of
Hamiltonian dynamics while allowing the broader contact-type Hamiltonian
motions.

 \noindent
\textbf{Contactization and Darboux Coordinates.}
A fundamental class of contact manifolds is obtained by contactizing the
cotangent bundle \(T^{*}Q\). More precisely, one considers the trivial real
line bundle over \(T^*Q\),
\begin{equation}\label{can-conison}
M=T^{*}Q\times \mathbb{R},
\end{equation}
with local coordinates \((q^{i},p_{i},z)\), where \((q^{i},p_{i})\) are the
canonical coordinates on \(T^{*}Q\) and \(z\) is the coordinate on the
\(\mathbb{R}\)-factor. On a local Darboux chart \(U_{\alpha}\), the standard
contact one-form is given by
\begin{equation}\label{eta-Q}
\eta\vert_{\alpha}
=
dz-\Theta_{Q}\vert_{\alpha}
=
dz-p_{i}dq^{i},
\qquad
\mathcal{R}^{\eta}\vert_{\alpha}
=
\frac{\partial}{\partial z},
\end{equation}
where \(\Theta_{Q}\) denotes, by abuse of notation, the pullback of the
Liouville one-form from \(T^{*}Q\) to \(T^{*}Q\times\mathbb{R}\). The associated
two-form is therefore
\begin{equation}\label{contact-two-form-local}
\Omega\vert_{\alpha}
=
d\eta\vert_{\alpha}
=
dq^{i}\wedge dp_{i}.
\end{equation}

Accordingly, the contact bivector field \(\Lambda^\eta\), introduced in
\eqref{Lambda-eta}, takes the local form
\begin{equation}\label{contact-bivector-local}
\Lambda^\eta\vert_{\alpha}
=
\frac{\partial}{\partial q^{i}}\wedge \frac{\partial}{\partial p_{i}}
+
p_{i}
\frac{\partial}{\partial z}\wedge \frac{\partial}{\partial p_{i}} .
\end{equation}
The contact bracket is then written as
\begin{equation}\label{contact-bracket-local} 
\{F,H\}^{\mathcal R} =
\frac{\partial F}{\partial q^{i}}
\frac{\partial H}{\partial p_{i}}
-
\frac{\partial F}{\partial p_{i}}
\frac{\partial H}{\partial q^{i}}
+
\left(
F-p_{i}\frac{\partial F}{\partial p_{i}}
\right)
\frac{\partial H}{\partial z}
-
\left(
H-p_{i}\frac{\partial H}{\partial p_{i}}
\right)
\frac{\partial F}{\partial z}. 
\end{equation}
Therefore, for a Hamiltonian function \(H\), the corresponding contact
Hamiltonian vector field is
\begin{equation}\label{contact-Hamiltonian-vector-field-local}
X_{H}^{\mathcal R}
=
\frac{\partial H}{\partial p_{i}}
\frac{\partial}{\partial q^{i}}
-
\left(
\frac{\partial H}{\partial q^{i}}
+
p_{i}\frac{\partial H}{\partial z}
\right)
\frac{\partial}{\partial p_{i}}
+
\left(
p_{i}\frac{\partial H}{\partial p_{i}}
-
H
\right)
\frac{\partial}{\partial z}.
\end{equation}
Hence, the contact Hamilton's equations are
\begin{equation}\label{contact-Hamilton-equations-local}
\frac{dq^{i}}{dt}
=
\frac{\partial H}{\partial p_{i}},
\qquad
\frac{dp_{i}}{dt}
=
-\frac{\partial H}{\partial q^{i}}
-
p_{i}\frac{\partial H}{\partial z},
\qquad
\frac{dz}{dt}
=
p_{i}\frac{\partial H}{\partial p_{i}}
-
H .
\end{equation}

\subsection{Conformal Analysis of Contact Manifolds.}

A \((2n+1)\)-dimensional manifold \(M\) is called a contact manifold if it
admits a smooth codimension one distribution, called a contact structure, which
is maximally non-integrable \cite{Arn,boothwang,grabowska2022contact,Takizawa}. Locally, this distribution may be
represented by one-forms \(\eta_{\alpha}\) on charts \(U_{\alpha}\) satisfying
\begin{equation}
(d\eta_{\alpha})^{n}\wedge \eta_{\alpha}\neq 0 .
\end{equation}
Thus, on each \(U_{\alpha}\), the local contact distribution is given by
\begin{equation}
\mathcal D|_{\alpha}
=
\ker \eta_{\alpha}
\subset TU_{\alpha}.
\end{equation}
In general, these local defining forms need not glue to a globally defined
one-form on \(M\). Hence the contact structure may exist globally as a
hyperplane distribution without being represented by a global contact form. A
contact manifold is called coorientable if there exists a global one-form
\(\eta\) such that
\begin{equation}
\mathcal D=\ker\eta,
\qquad
(d\eta)^n\wedge\eta\neq 0.
\end{equation}
The coorientable case has already been discussed in
Subsection~\ref{sec:Jacko}.

Motivated by the locally conformal cosymplectic construction above, we now introduce a locally conformal contact framework. The main point is to replace the
ordinary exterior derivative appearing in the coorientable contact case by the
Lichnerowicz--de~Rham differential. Thus, the contact-type two-form will not
be \(d\eta\), but rather its twisted analogue \(d_{\varphi}\eta\). This allows
one to keep track of the conformal deformation at the level of the associated
Hamiltonian-type geometry.

Let \(M\) be a \((2n+1)\)-dimensional manifold, and let \(\varphi\) be a closed
one-form on \(M\). A locally conformal contact structure, abbreviated as an
LCCon structure, is determined by 
the pair \((\Omega,\eta)\) satisfying the non-degeneracy condition
\begin{equation}\label{LCCon}
\Omega^{n}\wedge \eta \neq 0,
\qquad
\Omega=d_{\varphi}\eta .
\end{equation}
In this case, we write \((M,\Omega,\eta;\varphi)\) and call it a locally
conformal contact manifold. The closed one-form \(\varphi\) is called the Lee
form of the locally conformal contact structure.

Since
\begin{equation}
d_{\varphi}\eta=d\eta-\varphi\wedge\eta,
\end{equation}
the non-degeneracy condition in \eqref{LCCon} is equivalent to the usual
contact non-degeneracy condition
\begin{equation}
(d\eta)^n\wedge\eta\neq 0.
\end{equation}
Thus, the underlying hyperplane distribution \(\ker\eta\) is an ordinary
coorientable contact structure, while the locally conformal contact geometry is
encoded in the twisted two-form \(d_{\varphi}\eta\). In the particular case
\(\varphi=0\), the Lichnerowicz--de~Rham differential \(d_{\varphi}\) reduces
to the exterior derivative \(d\), and one recovers the coorientable contact
structure discussed in \eqref{cont-man-def}.

\noindent
\textbf{Where the Local Conformality Appears.}
We now explain how locally conformal contact structures \eqref{LCCon} arise from local
contact data. This point of view also clarifies the role of coorientability.
Recall that a contact structure is first of all a maximally non-integrable
hyperplane distribution. Locally, such a distribution is represented on each
chart \(U_{\alpha}\) by a contact form \(\eta_{\alpha}\), satisfying
\((d\eta_{\alpha})^{n}\wedge \eta_{\alpha}\neq 0\). Thus each
\((U_{\alpha},\eta_{\alpha})\) is a local coorientable contact manifold, and
\(\ker\eta_{\alpha}\subset TU_{\alpha}\) is a maximally non-integrable
codimension one distribution.

However, the existence of local contact forms does not by itself imply the
existence of a single global contact form. In other words, although the contact
geometry is available chart by chart, the local forms \(\eta_{\alpha}\) may fail
to glue directly to a globally defined one-form on \(M\). The locally conformal
contact construction organizes this failure by allowing the local contact forms
to agree only up to conformal rescaling.

Let \(M\) be a \((2n+1)\)-dimensional manifold with an open cover
\(\{U_\alpha\}\), and suppose that each chart \(U_\alpha\) is equipped with
a contact form \(\eta_\alpha\). We assume that there exist smooth functions
\(\sigma_\alpha\in C^\infty(U_\alpha)\) satisfying the compatibility
conditions
\begin{equation}\label{hoppa}
e^{\sigma_\alpha}\eta_\alpha
=
e^{\sigma_\beta}\eta_\beta,
\qquad
d\sigma_\alpha
=
d\sigma_\beta
\qquad
\text{on }
U_\alpha\cap U_\beta.
\end{equation}
We emphasize that the second equality in \eqref{hoppa} is an additional
compatibility condition; it does not follow from the conformal relation
between the contact forms alone. Equivalently, one may begin with a globally
defined closed Lee form \(\varphi\) and choose local potentials
\(\sigma_\alpha\) satisfying
\(\varphi|_\alpha=d\sigma_\alpha\).

The first compatibility condition in \eqref{hoppa} allows us to define a
global one-form \(\eta\) by
\begin{equation}\label{lcc-local-global-eta}
\eta|_\alpha
=
e^{\sigma_\alpha}\eta_\alpha.
\end{equation}
Indeed, these local representatives agree on every nonempty overlap.
Consequently, the hyperplane distribution
\[
\mathcal D=\ker\eta
\]
is globally cooriented.

The second compatibility condition in \eqref{hoppa} implies that the local
one-forms \(d\sigma_\alpha\) glue to a globally defined one-form
\(\varphi\), characterized by
\begin{equation}\label{lcc-local-global-lee}
\varphi|_\alpha
=
d\sigma_\alpha.
\end{equation}
Since
\[
d\varphi|_\alpha
=
d^2\sigma_\alpha
=
0
\]
on every chart, \(\varphi\) is closed. This globally defined closed one-form
is the Lee form of the LCCon structure.

It remains to verify that the local contact two-forms also glue with the
same conformal weight. Differentiating the first equality in
\eqref{hoppa}, we obtain
\begin{equation}
\begin{split}
d\left(e^{\sigma_\alpha}\eta_\alpha\right)
&=
e^{\sigma_\alpha}d\eta_\alpha
+
d\sigma_\alpha\wedge
e^{\sigma_\alpha}\eta_\alpha,
\\
d\left(e^{\sigma_\beta}\eta_\beta\right)
&=
e^{\sigma_\beta}d\eta_\beta
+
d\sigma_\beta\wedge
e^{\sigma_\beta}\eta_\beta.
\end{split}
\end{equation}
The left-hand sides are equal on \(U_\alpha\cap U_\beta\). Moreover, the
second terms on the right-hand sides agree by both compatibility conditions
in \eqref{hoppa}. Therefore,
\begin{equation}\label{LCCon-local-two-form-compatibility}
e^{\sigma_\alpha}d\eta_\alpha
=
e^{\sigma_\beta}d\eta_\beta
\qquad
\text{on }
U_\alpha\cap U_\beta.
\end{equation}
Hence, the local two-forms
\begin{equation}\label{lcc-local-global-omega}
\Omega|_\alpha
=
e^{\sigma_\alpha}d\eta_\alpha
\end{equation}
agree on chart overlaps and define a globally defined two-form
\(\Omega\).

Using
\(\eta|_\alpha=e^{\sigma_\alpha}\eta_\alpha\) and
\(\varphi|_\alpha=d\sigma_\alpha\), we obtain
\begin{equation} 
\Omega|_\alpha =
e^{\sigma_\alpha}d\eta_\alpha
 =
d\left(e^{\sigma_\alpha}\eta_\alpha\right)
-
d\sigma_\alpha\wedge
e^{\sigma_\alpha}\eta_\alpha
 =
d\eta|_\alpha
-
\varphi|_\alpha\wedge\eta|_\alpha. 
\end{equation}
Since this identity holds on every chart, it gives the global structural
relation
\begin{equation}\label{LCCon-global-twisted-relation}
\Omega
=
d\eta-\varphi\wedge\eta
=
d_\varphi\eta.
\end{equation}
This is precisely the structural relation required in the locally conformal
contact condition \eqref{LCCon}.

It remains to verify that the nondegeneracy condition is preserved under this
construction. We obtain
\begin{equation}
\Omega^{n}\vert_{\alpha}\wedge\eta\vert_{\alpha}
=
e^{(n+1)\sigma_{\alpha}}
(d\eta_{\alpha})^{n}\wedge\eta_{\alpha}.
\end{equation}
The conformal factor \(e^{(n+1)\sigma_{\alpha}}\) is nowhere vanishing, while
\((d\eta_{\alpha})^{n}\wedge\eta_{\alpha}\neq 0\) by the local contact
condition. Therefore
\begin{equation}
\Omega^{n}\wedge\eta\neq 0
\end{equation}
globally. Consequently, the pair \((\Omega,\eta)\), together with the closed Lee
form \(\varphi\), satisfies exactly the defining conditions of a locally
conformal contact structure in \eqref{LCCon}.

Consequently, the local contact data determine a globally defined locally
conformally contact structure \((M,\Omega,\eta;\varphi)\), with
\(\Omega=d_{\varphi}\eta\) and \(d\varphi=0\). The terminology ``locally
conformally contact'' is therefore literal: after suitable local conformal
rescalings, the structure is represented on each chart by an ordinary contact
form. Global coorientability follows from the assumed compatibility of the weighted local contact forms, which produces a global nowhere-vanishing one-form \(\eta\).  

\noindent
\textbf{From Local Jacobi Structures to a Global Jacobi Structure.}
We now pass from the Jacobi structures associated with the local contact
data to a globally defined Jacobi structure. Each local contact manifold
\((U_\alpha,\eta_\alpha)\) determines a Jacobi pair
\[
\left(
\Lambda_\alpha^\eta,
\mathcal R_\alpha^\eta
\right),
\]
where \(\Lambda_\alpha^\eta\) is the local contact bivector and
\(\mathcal R_\alpha^\eta\) is the local contact Reeb vector field. With the
sign convention adopted in this paper, the local Jacobi identities are
\begin{equation}\label{local-LCCon-Jacobi-identities}
[\Lambda_\alpha^\eta,\Lambda_\alpha^\eta]_{SN}
=
-2\mathcal R_\alpha^\eta\wedge\Lambda_\alpha^\eta,
\qquad
[\mathcal R_\alpha^\eta,\Lambda_\alpha^\eta]_{SN}
=
0.
\end{equation}

Under the equal-weight conformal globalization, the local bivector and Reeb
vector fields transform according to
\begin{equation}\label{LCCon-local-bivector-Reeb}
\Lambda^\eta|_\alpha
=
e^{-\sigma_\alpha}\Lambda_\alpha^\eta,
\qquad
\mathcal R^\eta|_\alpha
=
e^{-\sigma_\alpha}\mathcal R_\alpha^\eta.
\end{equation}
Equivalently,
\[
\Lambda_\alpha^\eta
=
e^{\sigma_\alpha}\Lambda^\eta|_\alpha,
\qquad
\mathcal R_\alpha^\eta
=
e^{\sigma_\alpha}\mathcal R^\eta|_\alpha.
\]
Although the local pairs are Jacobi, the conformal factors are generally
nonconstant. Their derivatives therefore contribute additional terms to the
global Schouten--Nijenhuis identities.

We begin with the second local Jacobi identity in
\eqref{local-LCCon-Jacobi-identities}. Using
\eqref{LCCon-local-bivector-Reeb} and the Schouten--Nijenhuis product
identities, we obtain
\begin{equation}
\begin{split}
0
&=
[\mathcal R_\alpha^\eta,\Lambda_\alpha^\eta]_{SN}
\\
&=
\left[
e^{\sigma_\alpha}\mathcal R^\eta|_\alpha,
e^{\sigma_\alpha}\Lambda^\eta|_\alpha
\right]_{SN}
\\
&=
e^{2\sigma_\alpha}
\left(
\iota_{d\sigma_\alpha}
\bigl(
\Lambda^\eta|_\alpha\wedge\mathcal R^\eta|_\alpha
\bigr)
-
[\Lambda^\eta|_\alpha,\mathcal R^\eta|_\alpha]_{SN}
\right).
\end{split}
\end{equation}
Consequently,
\[
[\Lambda^\eta|_\alpha,\mathcal R^\eta|_\alpha]_{SN}
=
\iota_{d\sigma_\alpha}
\bigl(
\Lambda^\eta|_\alpha\wedge\mathcal R^\eta|_\alpha
\bigr).
\]
Since \(d\sigma_\alpha=\varphi|_\alpha\), this identity globalizes as
\begin{equation}\label{LCCon-Jacobi-second-pre}
[\Lambda^\eta,\mathcal R^\eta]_{SN}
=
\iota_\varphi
\bigl(
\Lambda^\eta\wedge\mathcal R^\eta
\bigr).
\end{equation}
Thus, the naively globalized pair
\((\Lambda^\eta,\mathcal R^\eta)\) does not, in general, satisfy the second
Jacobi identity. Its obstruction is determined precisely by the Lee form.

We next consider the first local Jacobi identity. From
\eqref{LCCon-local-bivector-Reeb}, the conformal transformation rule for the
Schouten--Nijenhuis bracket gives
\begin{equation}
\begin{split}
[\Lambda_\alpha^\eta,\Lambda_\alpha^\eta]_{SN}
&=
\left[
e^{\sigma_\alpha}\Lambda^\eta|_\alpha,
e^{\sigma_\alpha}\Lambda^\eta|_\alpha
\right]_{SN}
\\
&=
e^{2\sigma_\alpha}
\left(
[\Lambda^\eta|_\alpha,\Lambda^\eta|_\alpha]_{SN}
+
2\,
\iota_{d\sigma_\alpha}\Lambda^\eta|_\alpha
\wedge\Lambda^\eta|_\alpha
\right).
\end{split}
\end{equation}
On the other hand, the first identity in
\eqref{local-LCCon-Jacobi-identities} gives
\begin{equation}
[\Lambda_\alpha^\eta,\Lambda_\alpha^\eta]_{SN}
=
-2e^{2\sigma_\alpha}
\mathcal R^\eta|_\alpha\wedge\Lambda^\eta|_\alpha.
\end{equation}
Comparing these expressions yields
\begin{equation}
\begin{split}
[\Lambda^\eta|_\alpha,\Lambda^\eta|_\alpha]_{SN}
&=
-2
\left(
\iota_{d\sigma_\alpha}\Lambda^\eta|_\alpha
+
\mathcal R^\eta|_\alpha
\right)
\wedge\Lambda^\eta|_\alpha.
\end{split}
\end{equation}
Using \(d\sigma_\alpha=\varphi|_\alpha\), we obtain the global relation
\begin{equation}\label{LCCon-Jacobi-first-global}
[\Lambda^\eta,\Lambda^\eta]_{SN}
=
-2
\left(
\iota_\varphi\Lambda^\eta
+
\mathcal R^\eta
\right)
\wedge\Lambda^\eta.
\end{equation}

Equation \eqref{LCCon-Jacobi-first-global} identifies the correction required
by the conformal globalization. We therefore introduce the modified
Reeb-type vector field
\begin{equation}\label{LCCon-modified-Reeb}
\widetilde{\mathcal R}^\eta
=
\iota_\varphi\Lambda^\eta
+
\mathcal R^\eta.
\end{equation}
In terms of this vector field,
\eqref{LCCon-Jacobi-first-global} becomes
\begin{equation}\label{LCCon-Jacobi-first-modified}
[\Lambda^\eta,\Lambda^\eta]_{SN}
=
-2\widetilde{\mathcal R}^\eta\wedge\Lambda^\eta.
\end{equation}

It remains to verify the second Jacobi identity for the modified pair. Since
the Lee form is closed, it is locally exact. The contraction identity for
the Schouten--Nijenhuis bracket therefore gives
\begin{equation}\label{LCCon-contraction-SN-identity}
[\Lambda^\eta,\iota_\varphi\Lambda^\eta]_{SN}
=
\frac{1}{2}
\iota_\varphi
[\Lambda^\eta,\Lambda^\eta]_{SN}.
\end{equation}
Using \eqref{LCCon-Jacobi-second-pre},
\eqref{LCCon-Jacobi-first-global}, and
\eqref{LCCon-contraction-SN-identity}, we compute
\begin{equation}
\begin{split}
[\Lambda^\eta,\widetilde{\mathcal R}^\eta]_{SN}
&=
[\Lambda^\eta,\iota_\varphi\Lambda^\eta]_{SN}
+
[\Lambda^\eta,\mathcal R^\eta]_{SN}
\\
&=
\frac{1}{2}
\iota_\varphi
[\Lambda^\eta,\Lambda^\eta]_{SN}
+
\iota_\varphi
\bigl(
\Lambda^\eta\wedge\mathcal R^\eta
\bigr)
\\
&=
-
\iota_\varphi
\left(
\left(
\iota_\varphi\Lambda^\eta+\mathcal R^\eta
\right)
\wedge\Lambda^\eta
\right)
+
\iota_\varphi
\bigl(
\Lambda^\eta\wedge\mathcal R^\eta
\bigr)
\\
&=
-
\iota_\varphi
\bigl(
\iota_\varphi\Lambda^\eta\wedge\Lambda^\eta
\bigr)
-
\iota_\varphi
\bigl(
\mathcal R^\eta\wedge\Lambda^\eta
\bigr)
+
\iota_\varphi
\bigl(
\Lambda^\eta\wedge\mathcal R^\eta
\bigr)
\\
&=
0.
\end{split}
\end{equation}
Indeed,
\[
\iota_\varphi
\bigl(
\iota_\varphi\Lambda^\eta\wedge\Lambda^\eta
\bigr)
=
0,
\]
whereas
\[
\mathcal R^\eta\wedge\Lambda^\eta
=
\Lambda^\eta\wedge\mathcal R^\eta,
\]
so the remaining contraction terms cancel. Equivalently,
\begin{equation}\label{LCCon-Jacobi-second-modified}
[\widetilde{\mathcal R}^\eta,\Lambda^\eta]_{SN}
=
0.
\end{equation}

The preceding calculation establishes the global Jacobi structure of the
LCCon geometry.

\begin{proposition}
\label{prop:LCCon-global-Jacobi}
Let \((M,\Omega,\eta;\varphi)\) be an LCCon manifold obtained by the
equal-weight conformal globalization of local contact structures. Then
\begin{equation}
\big(
\Lambda^\eta,
\widetilde{\mathcal R}^\eta
\big),
\qquad
\widetilde{\mathcal R}^\eta
=
\iota_\varphi\Lambda^\eta+\mathcal R^\eta,
\end{equation}
defines a global Jacobi structure on \(M\). More precisely,
\begin{equation}\label{LCCon-global-Jacobi-identities}
[\Lambda^\eta,\Lambda^\eta]_{SN}
=
-2\widetilde{\mathcal R}^\eta\wedge\Lambda^\eta,
\qquad
[\widetilde{\mathcal R}^\eta,\Lambda^\eta]_{SN}
=
0.
\end{equation}
\end{proposition}

Thus, the Lee form measures the failure of the naively globalized pair
\((\Lambda^\eta,\mathcal R^\eta)\) to satisfy the Jacobi identities, whereas
the modified vector field \(\widetilde{\mathcal R}^\eta\) absorbs this
failure. In particular, when \(\varphi=0\), the correction term vanishes and
the ordinary contact Jacobi pair
\((\Lambda^\eta,\mathcal R^\eta)\) is recovered.

\begin{proposition}
\label{prop:LCCon-underlying-contact}
Let \((M,\Omega,\eta;\varphi)\) be an LCCon manifold, with
\[
\Omega=d_\varphi\eta=d\eta-\varphi\wedge\eta.
\]
Then \(\eta\) is an ordinary contact form. Moreover, the bivector field
\(\Lambda^\eta\) associated with the odd pair \((\Omega,\eta)\) coincides
with the contact bivector determined by \((d\eta,\eta)\), and the modified
vector field
\[
\widetilde{\mathcal R}^\eta
=
\iota_\varphi\Lambda^\eta+\mathcal R^\eta
\]
is precisely the contact Reeb vector field of \(\eta\).
\end{proposition}

The proposition shows that the LCCon construction does not replace the
underlying contact geometry by a different geometric category. Rather, it
provides a conformally enriched odd-pair description of that geometry. The
twisted two-form \(\Omega=d_\varphi\eta\) and the odd-pair Reeb vector field
\(\mathcal R^\eta\) encode the conformal globalization, while
\((\Lambda^\eta,\widetilde{\mathcal R}^\eta)\) recovers the canonical Jacobi
structure of the global contact form \(\eta\).

\subsection{Conformal Analysis of Contact Dynamics}

The distinction between the ordinary contact formulation and its locally
conformal description becomes especially visible at the level of Hamiltonian
dynamics. On each chart \(U_\alpha\), the dynamics is governed by an ordinary
contact Hamiltonian system associated with the local contact form
\(\eta_\alpha\). Globally, the conformal gluing of these local systems is
encoded by the Lee form \(\varphi\) and the Lichnerowicz--de~Rham differential
\(d_\varphi\). In particular, this globalization distinguishes the Reeb-type
vector field associated with the non-degenerate odd pair \((\Omega,\eta)\)
from the Reeb vector field governing the resulting Jacobi Hamiltonian
dynamics.

As in the LCCos case, there are two natural approaches to LCCon Hamiltonian
dynamics. The first directly replaces the exterior derivative \(d\) in the
ordinary contact Hamiltonian equations by the Lichnerowicz--de~Rham
differential \(d_\varphi\). The second starts with ordinary contact
Hamiltonian systems on local charts and globalizes them through the conformal
transition rules. We now present these constructions separately.

\noindent
\textbf{First Approach: Direct Lichnerowicz--de~Rham Deformation.}
Let \((M,\Omega,\eta;\varphi)\) be an LCCon manifold, where
\[
\Omega=d_\varphi\eta,
\qquad
d\varphi=0,
\]
and let \(H\in C^\infty(M)\). Replacing \(dH\) in the ordinary contact
Hamiltonian equations by \(d_\varphi H\), we define the LCCon Hamiltonian
vector field \(X_H^{\mathcal R}\) through
\begin{equation}\label{LCCon-Ham-VF}
\iota_{X_H^{\mathcal R}}\Omega
=
d_\varphi H
-
\left\langle
d_\varphi H,\mathcal R^\eta
\right\rangle\eta,
\qquad
\iota_{X_H^{\mathcal R}}\eta=-H,
\end{equation}
where \(\mathcal R^\eta\) is the Reeb-type vector field associated with the
non-degenerate odd pair \((\Omega,\eta)\), characterized by
\[
\iota_{\mathcal R^\eta}\Omega=0,
\qquad
\iota_{\mathcal R^\eta}\eta=1.
\]

Equivalently, using the bivector field \(\Lambda^\eta\), we obtain
\begin{equation}\label{lcc-global-ham-vector}
X_H^{\mathcal R}
=
\sharp_{\Lambda^\eta}(d_\varphi H)
-
H\mathcal R^\eta.
\end{equation}
Expanding \(d_\varphi H=dH-H\varphi\) gives
\begin{equation}\label{LCCon-Ham-modified-Reeb-form}
\begin{split}
X_H^{\mathcal R}
&=
\sharp_{\Lambda^\eta}(dH-H\varphi)
-
H\mathcal R^\eta
\\
&=
\sharp_{\Lambda^\eta}(dH)
-
H\left(
\iota_\varphi\Lambda^\eta+\mathcal R^\eta
\right).
\end{split}
\end{equation}
Recalling the modified Reeb-type vector field
\(\widetilde{\mathcal R}^\eta\) introduced in
\eqref{LCCon-modified-Reeb}, we arrive at
\begin{equation}\label{LCCon-Ham-Jacobi-form}
X_H^{\mathcal R}
=
\sharp_{\Lambda^\eta}(dH)
-
H\widetilde{\mathcal R}^\eta.
\end{equation}
By Proposition~\ref{prop:LCCon-global-Jacobi}, this is precisely the
Hamiltonian vector field associated with the global Jacobi structure
\((\Lambda^\eta,\widetilde{\mathcal R}^\eta)\).

\noindent
\textbf{Second Approach: Globalization of Local Contact Hamiltonian
Dynamics.}
Let \(H_\alpha\) be a local Hamiltonian function on the contact chart
\((U_\alpha,\eta_\alpha)\), and denote the corresponding local contact
Hamiltonian vector field by \(X_{H_\alpha}^{\mathcal R}\). It is determined
by
\begin{equation}\label{local-contact-ham}
\iota_{X_{H_\alpha}^{\mathcal R}}d\eta_\alpha
=
dH_\alpha
-
\mathcal R_\alpha^\eta(H_\alpha)\eta_\alpha,
\qquad
\iota_{X_{H_\alpha}^{\mathcal R}}\eta_\alpha
=
-H_\alpha,
\end{equation}
where \(\mathcal R_\alpha^\eta\) is the contact Reeb vector field associated
with \(\eta_\alpha\).

We assume that the local Hamiltonians transform with the same conformal
weight as the local geometric data:
\begin{equation}\label{lcc-local-comp-h}
e^{\sigma_\alpha}H_\alpha
=
e^{\sigma_\beta}H_\beta
\qquad
\text{on }
U_\alpha\cap U_\beta.
\end{equation}
They therefore determine a globally defined Hamiltonian function
\(H\in C^\infty(M)\), locally represented by
\begin{equation}\label{lcc-local-global-h}
H|_\alpha
=
e^{\sigma_\alpha}H_\alpha.
\end{equation}
Using \(\varphi|_\alpha=d\sigma_\alpha\), we obtain
\begin{equation}\label{LCCon-twisted-H-local}
\left.d_\varphi H\right|_\alpha
=
d\left(e^{\sigma_\alpha}H_\alpha\right)
-
d\sigma_\alpha\,e^{\sigma_\alpha}H_\alpha
=
e^{\sigma_\alpha}dH_\alpha.
\end{equation}

The global geometric data are locally represented by
\begin{equation}
\eta|_\alpha
=
e^{\sigma_\alpha}\eta_\alpha,
\qquad
\Omega|_\alpha
=
e^{\sigma_\alpha}d\eta_\alpha,
\qquad
\mathcal R^\eta|_\alpha
=
e^{-\sigma_\alpha}\mathcal R_\alpha^\eta.
\end{equation}
Here \(\mathcal R^\eta\) is the Reeb-type vector field of the global odd pair
\((\Omega,\eta)\). It should not be confused with the modified vector field
\(\widetilde{\mathcal R}^\eta\) entering the global Jacobi structure.

We first examine the equation involving the contact form. Using
\eqref{local-contact-ham}, we obtain
\begin{equation}
\iota_{X_{H_\alpha}^{\mathcal R}}\eta|_\alpha
=
e^{\sigma_\alpha}
\iota_{X_{H_\alpha}^{\mathcal R}}\eta_\alpha
=
-e^{\sigma_\alpha}H_\alpha
=
-H|_\alpha.
\end{equation}
Thus, the local normalization equations produce the global identity
\begin{equation}\label{lcc-global-ham-second}
\iota_{X_H^{\mathcal R}}\eta=-H.
\end{equation}

For the two-form equation, we compute
\begin{equation}
\iota_{X_{H_\alpha}^{\mathcal R}}\Omega|_\alpha
=
e^{\sigma_\alpha}
\iota_{X_{H_\alpha}^{\mathcal R}}d\eta_\alpha
=
e^{\sigma_\alpha}
\left(
dH_\alpha
-
\mathcal R_\alpha^\eta(H_\alpha)\eta_\alpha
\right).
\end{equation}
The first term satisfies
\[
e^{\sigma_\alpha}dH_\alpha
=
\left.d_\varphi H\right|_\alpha.
\]
Moreover,
\begin{equation} 
\left.
\left\langle
d_\varphi H,\mathcal R^\eta
\right\rangle
\right|_\alpha
=
\left\langle
e^{\sigma_\alpha}dH_\alpha,
e^{-\sigma_\alpha}\mathcal R_\alpha^\eta
\right\rangle=
\mathcal R_\alpha^\eta(H_\alpha). 
\end{equation}
Consequently,
\begin{equation}
\iota_{X_{H_\alpha}^{\mathcal R}}\Omega|_\alpha
=
\left.
\left(
d_\varphi H
-
\left\langle
d_\varphi H,\mathcal R^\eta
\right\rangle\eta
\right)
\right|_\alpha.
\end{equation}
Therefore, every local contact Hamiltonian vector field satisfies the
restriction of the global equation
\begin{equation}\label{lcc-global-ham-first}
\iota_{X_H^{\mathcal R}}\Omega
=
d_\varphi H
-
\left\langle
d_\varphi H,\mathcal R^\eta
\right\rangle\eta.
\end{equation}
By the non-degeneracy of the odd pair \((\Omega,\eta)\), the system formed by
\eqref{lcc-global-ham-first} and \eqref{lcc-global-ham-second} has a unique
solution. Hence, the local contact Hamiltonian vector fields agree on chart
overlaps and define a global vector field.

The equivalence of the two constructions is summarized in the following
result.

\begin{proposition}
\label{prop:LCCon-Hamiltonian-equivalence}
Let \((M,\Omega,\eta;\varphi)\) be an LCCon manifold obtained by conformally
gluing local contact structures, and suppose that the local Hamiltonians
\(\{H_\alpha\}\) satisfy the compatibility condition
\eqref{lcc-local-comp-h}. Then the corresponding local contact Hamiltonian
vector fields glue directly to a global vector field:
\begin{equation}
X_H^{\mathcal R}|_\alpha
=
X_{H_\alpha}^{\mathcal R}.
\end{equation}
This global vector field coincides with the vector field obtained by direct
Lichnerowicz--de~Rham deformation and is characterized by
\begin{equation}
\iota_{X_H^{\mathcal R}}\Omega
=
d_\varphi H
-
\left\langle
d_\varphi H,\mathcal R^\eta
\right\rangle\eta,
\qquad
\iota_{X_H^{\mathcal R}}\eta=-H.
\end{equation}
Equivalently,
\begin{equation}
X_H^{\mathcal R}
=
\sharp_{\Lambda^\eta}(dH)
-
H\widetilde{\mathcal R}^\eta,
\qquad
\widetilde{\mathcal R}^\eta
=
\iota_\varphi\Lambda^\eta+\mathcal R^\eta.
\end{equation}
Thus, direct deformation and conformal globalization provide two equivalent
formulations of the same LCCon Hamiltonian dynamics.
\end{proposition}

The proposition shows that the Lee-form correction is not introduced
externally into the local contact dynamics. Rather, it arises intrinsically
from the conformal globalization process and appears globally through the
modification of the odd-pair Reeb vector field. When \(\varphi=0\), one has
\(\widetilde{\mathcal R}^\eta=\mathcal R^\eta\), and both constructions
reduce to ordinary contact Hamiltonian dynamics.

\noindent
\textbf{Coordinate Formulas.}
We now derive the coordinate expression of the locally conformally contact
Hamiltonian dynamics. The aim is to display explicitly how the Lee form
modifies the ordinary contact Hamilton's equations.

Let \((q^i,p_i,z)\) be local coordinates on a chart of
\(T^*Q\times\mathbb R\). We write the closed Lee one-form as
\begin{equation}
\varphi=\varphi_i\,dq^i+\psi\,dz.
\end{equation}
Here \(\varphi_i\) are the components of the Lee form along the base
directions, whereas \(\psi\) is its component along the contact coordinate
\(z\).

Let \((q^i,r_i,\zeta)\) be local contact Darboux coordinates on \(U_\alpha\),
so that
\begin{equation}
\eta_\alpha=d\zeta-r_i\,dq^i.
\end{equation}
We assume that the local contact coordinates are related to the global
coordinates by
\begin{equation}
\zeta=e^{-\sigma_\alpha}z,
\qquad
r_i=e^{-\sigma_\alpha}p_i.
\end{equation}
Thus, the global one-form is obtained from
\(\eta=e^{\sigma_\alpha}\eta_\alpha\). Since
\(\varphi|_\alpha=d\sigma_\alpha\), a direct computation gives
\begin{equation}\label{LCCon-eta-coord}
\eta
=
(1-z\psi)\,dz
-
(p_i+z\varphi_i)\,dq^i.
\end{equation}
This formula exhibits the first effect of the Lee form: the coefficient of
\(dz\) is modified by \(\psi\), whereas the coefficient of \(dq^i\) is
shifted from \(p_i\) to \(p_i+z\varphi_i\).

Equivalently, the same one-form can be written intrinsically as
\begin{equation}\label{LCCon-eta-intrinsic}
\eta=d_\varphi z-\Theta_Q,
\end{equation}
where \(\Theta_Q=p_i\,dq^i\) denotes the pull-back of the Liouville one-form.
Indeed, \(d_\varphi z=dz-z\varphi\), and substitution of
\(\varphi=\varphi_i\,dq^i+\psi\,dz\) gives
\eqref{LCCon-eta-coord}.

Applying the twisted differential once more, we obtain
\begin{equation}\label{LCCon-Omega-coord}
\Omega
=
dq^i\wedge dp_i
-
\frac{1}{2}C_{ij}\,dq^i\wedge dq^j
+
p_i\psi\,dz\wedge dq^i,
\end{equation}
where
\begin{equation}
C_{ij}
=
\varphi_jp_i-\varphi_ip_j.
\end{equation}
Thus, \(\Omega\) consists of the canonical cotangent contribution together
with two Lee-form corrections. The first is the skew term involving
\(C_{ij}\), generated by the base components \(\varphi_i\), whereas the
second is the mixed term \(p_i\psi\,dz\wedge dq^i\), which is generated by
the \(dz\)-component of the Lee form.

\noindent
\textbf{The Odd-Pair and Jacobi Reeb Vector Fields.}
At this point, two related but distinct Reeb-type vector fields must be
distinguished. The non-degenerate odd pair \((\Omega,\eta)\) determines
\(\mathcal R^\eta\) through
\begin{equation}
\iota_{\mathcal R^\eta}\Omega=0,
\qquad
\iota_{\mathcal R^\eta}\eta=1.
\end{equation}
Solving these equations gives, on the domain where \(1-z\psi\neq0\),
\begin{equation}\label{LCCon-Reeb-coord}
\mathcal R^\eta
=
\frac{p_i\psi}{1-z\psi}
\frac{\partial}{\partial p_i}
+
\frac{1}{1-z\psi}
\frac{\partial}{\partial z}.
\end{equation}

The global Jacobi structure, however, is determined by
\begin{equation}
\widetilde{\mathcal R}^\eta
=
\iota_\varphi\Lambda^\eta+\mathcal R^\eta.
\end{equation}
This vector field is, in fact, the contact Reeb vector field of the global
contact form \(\eta\). Indeed, since
\[
\Omega=d_\varphi\eta=d\eta-\varphi\wedge\eta,
\]
we have
\[
d\eta=\Omega+\varphi\wedge\eta.
\]
Using
\[
\iota_{\sharp_{\Lambda^\eta}(\varphi)}\Omega
=
\varphi-\varphi(\mathcal R^\eta)\eta,
\qquad
\eta\bigl(\sharp_{\Lambda^\eta}(\varphi)\bigr)=0,
\]
one obtains
\begin{equation}
\iota_{\widetilde{\mathcal R}^\eta}\eta=1,
\qquad
\iota_{\widetilde{\mathcal R}^\eta}d\eta=0.
\end{equation}
Consequently, \(\mathcal R^\eta\) is the Reeb-type vector field of the odd
pair \((\Omega,\eta)\), whereas \(\widetilde{\mathcal R}^\eta\) is the Reeb
vector field entering the global contact Jacobi structure.

Its coordinate expression is
\begin{equation}\label{LCCon-modified-Reeb-coord}
\widetilde{\mathcal R}^\eta
=
-\frac{\varphi_i}{1-z\psi}
\frac{\partial}{\partial p_i}
+
\frac{1}{1-z\psi}
\frac{\partial}{\partial z}.
\end{equation}

We now determine the Hamiltonian vector field. For
\(H\in C^\infty(M)\), the LCCon Hamiltonian vector field is defined by
\begin{equation}\label{LCCon-Ham-def-coord}
\iota_{X_H^{\mathcal R}}\Omega
=
d_\varphi H
-
\left\langle
d_\varphi H,\mathcal R^\eta
\right\rangle\eta,
\qquad
\iota_{X_H^{\mathcal R}}\eta=-H.
\end{equation}
Substituting the coordinate expressions of
\(\eta\), \(\Omega\), and \(\mathcal R^\eta\) gives
\begin{equation}\label{LCCon-Ham-VF-coord}
\begin{split}
X_H^{\mathcal R}
&=
\frac{\partial H}{\partial p_i}
\frac{\partial}{\partial q^i}
+
\left(
-\frac{\partial H}{\partial q^i}
+
\frac{H\varphi_i}{1-z\psi}
-
\frac{p_i+z\varphi_i}{1-z\psi}
\frac{\partial H}{\partial z}
+
\frac{C_{ij}}{1-z\psi}
\frac{\partial H}{\partial p_j}
\right)
\frac{\partial}{\partial p_i}
\\
&\quad
+
\left(
\frac{p_i+z\varphi_i}{1-z\psi}
\frac{\partial H}{\partial p_i}
-
\frac{H}{1-z\psi}
\right)
\frac{\partial}{\partial z}.
\end{split}
\end{equation}
As a consistency check, setting \(H=1\) in this expression gives
\begin{equation}
X_1^{\mathcal R}
=
-\widetilde{\mathcal R}^\eta,
\end{equation}
as expected from the Jacobi representation
\eqref{LCCon-Ham-Jacobi-form}.

The $q^i$-equation retains its familiar form $\dot q^i = \partial H/\partial p_i$. The Lee-form modifications appear in
the \(p_i\)- and \(z\)-components through the drift term \(H\varphi_i\), the
shifted coefficient \(p_i+z\varphi_i\), the skew term \(C_{ij}\), and the
common denominator \(1-z\psi\).

Therefore, the corresponding LCCon Hamilton's equations
are
\begin{equation}\label{LCCon-Ham-eq-coord}
\begin{split}
\frac{dq^i}{dt}
&=
\frac{\partial H}{\partial p_i},
\\
\frac{dp_i}{dt}
&=
-\frac{\partial H}{\partial q^i}
+
\frac{H\varphi_i}{1-z\psi}
-
\frac{p_i+z\varphi_i}{1-z\psi}
\frac{\partial H}{\partial z}
+
\frac{C_{ij}}{1-z\psi}
\frac{\partial H}{\partial p_j},
\\
\frac{dz}{dt}
&=
\frac{p_i+z\varphi_i}{1-z\psi}
\frac{\partial H}{\partial p_i}
-
\frac{H}{1-z\psi}.
\end{split}
\end{equation}
When the Lee form vanishes, that is,
\(\varphi_i=0\) and \(\psi=0\), all the conformal correction terms disappear,
the two Reeb fields coincide, and the standard contact Hamilton equations are
recovered.

A particularly important simplification occurs when the Lee form has no
component along the contact coordinate, namely when \(\psi=0\). In this case,
\[
\varphi=\varphi_i\,dq^i,
\qquad
\eta=dz-(p_i+z\varphi_i)dq^i.
\]
The odd-pair Reeb vector field becomes purely vertical:
\begin{equation}
\mathcal R^\eta
=
\frac{\partial}{\partial z}.
\end{equation}
Nevertheless, the Jacobi Reeb vector field remains conformally modified:
\begin{equation}
\widetilde{\mathcal R}^\eta
=
\frac{\partial}{\partial z}
-
\varphi_i\frac{\partial}{\partial p_i}.
\end{equation}
Thus, although the odd-pair Reeb field reduces to its ordinary vertical form,
the Reeb field governing the Jacobi Hamiltonian dynamics retains a
Lee-form-dependent momentum component.

In this case, the Hamiltonian vector field becomes
\begin{equation}\label{LCCon-Ham-VF-psi-zero}
\begin{split}
X_H^{\mathcal R}
&=
\frac{\partial H}{\partial p_i}
\frac{\partial}{\partial q^i} 
+
\left(
-\frac{\partial H}{\partial q^i}
+
H\varphi_i
-
(p_i+z\varphi_i)
\frac{\partial H}{\partial z}
+
C_{ij}
\frac{\partial H}{\partial p_j}
\right)
\frac{\partial}{\partial p_i}
\\
&\quad
+
\left(
(p_i+z\varphi_i)
\frac{\partial H}{\partial p_i}
-
H
\right)
\frac{\partial}{\partial z}.
\end{split}
\end{equation}
Accordingly, Hamilton's equations reduce to
\begin{equation}\label{LCCon-Ham-eq-psi-zero}
\begin{split}
\frac{dq^i}{dt}
&=
\frac{\partial H}{\partial p_i},
\\
\frac{dp_i}{dt}
&=
-\frac{\partial H}{\partial q^i}
+
H\varphi_i
-
(p_i+z\varphi_i)
\frac{\partial H}{\partial z}
+
C_{ij}
\frac{\partial H}{\partial p_j},
\\
\frac{dz}{dt}
&=
(p_i+z\varphi_i)
\frac{\partial H}{\partial p_i}
-
H.
\end{split}
\end{equation}
Therefore, even in the case \(\psi=0\), the LCCon representation does not
reduce to the undeformed coordinate description. The Lee form continues to
modify the equations through the drift \(H\varphi_i\), the shifted contact
momentum \(p_i+z\varphi_i\), the skew response \(C_{ij}\), and the momentum
component of the Jacobi Reeb vector field. 
 
\subsection{Dissipative Particle Dynamics in a LC Contact Environment}\label{sec:LCCon-particle-example}

In this subsection, we present a mechanical model illustrating how the LCCon
framework modifies the standard contact description of dissipative particle
dynamics. The model is based on a contact Hamiltonian with a Riemannian
kinetic energy, a potential energy, and a linear dependence on the action
variable. In the absence of the Lee form, the resulting equations reduce to
the usual contact Hamiltonian description of a linearly damped mechanical
system. When the Lee form is nonzero, the ordinary dissipative dynamics is
accompanied by additional locally conformal force terms.

Let \(Q\) be a three-dimensional manifold with local coordinates
\(\mathbf q=(q^1,q^2,q^3)\), and consider the locally conformally contact
phase space \(T^*Q\times\mathbb R\), with induced local coordinates
\((\mathbf q,\mathbf p,z)\). Let
\(\mathcal G(\mathbf q)=[\mathcal G_{ij}(\mathbf q)]\) be a Riemannian metric
matrix on \(Q\), and let
\[
G(\mathbf q)=\mathcal G(\mathbf q)^{-1}
=
[G^{ij}(\mathbf q)]
\]
be its inverse matrix. Note that the metric is assumed to be independent of
time.

We write the Lee one-form as
\begin{equation}\label{LCCon-Lee-form-Riemannian}
\varphi
=
\boldsymbol{\varphi}(\mathbf q,z)\cdot d\mathbf q
+
\psi(\mathbf q,z)dz,
\end{equation}
where
\(\boldsymbol{\varphi}=(\varphi_1,\varphi_2,\varphi_3)\).
On the open region where \(1-z\psi\neq 0\), consider the Hamiltonian
\begin{equation}\label{LCCon-dissipative-particle-H-Riemannian}
H(\mathbf q,\mathbf p,z)
=
\frac{1}{2}\mathbf p\cdot G(\mathbf q)\mathbf p
+
V(\mathbf q)
+
\gamma z .
\end{equation}
We assume \(\gamma>0\). The term \(\gamma z\) produces the
velocity-linear contribution with coefficient
\(\gamma/(1-z\psi)\); this contribution acts as damping on
the region \(1-z\psi>0\).

Substituting \eqref{LCCon-dissipative-particle-H-Riemannian} into the general
LCCon Hamilton's equations \eqref{LCCon-Ham-eq-coord}, we obtain
\begin{equation}\label{LCCon-dissipative-particle-Ham-eq-Riemannian}
\begin{split}
\dot{\mathbf q}
&=
G\mathbf p,
\\
\dot{\mathbf p}
&=
-
\nabla_{\mathbf q}H
+
\frac{1}{1-z\psi}
\left(
(\boldsymbol{\varphi}\times\mathbf p)\times G\mathbf p
+
\boldsymbol{\varphi}
\left(
\frac{1}{2}\mathbf p\cdot G\mathbf p
+
V(\mathbf q)
\right)
-
\gamma\mathbf p
\right),
\\
\dot z
&=
\frac{1}{1-z\psi}
\left(
\frac{1}{2}\mathbf p\cdot G\mathbf p
-
V(\mathbf q)
-
\gamma z
+
z\,\boldsymbol{\varphi}\cdot G\mathbf p
\right).
\end{split}
\end{equation}
Here
\begin{equation}\label{LCCon-gradient-H-Riemannian}
\nabla_{\mathbf q}H
=
\frac{1}{2}
\nabla_{\mathbf q}
\bigl(\mathbf p\cdot G(\mathbf q)\mathbf p\bigr)
+
\nabla_{\mathbf q}V(\mathbf q).
\end{equation}
The term 
$(\boldsymbol{\varphi}\times\mathbf p)\times G\mathbf p$ 
is precisely the LCCon counterpart of the skew contribution already appearing
in the LCCos Riemannian example. Hence we do not repeat its derivation here.

We consider the following identifications between the generalized velocities
and the momenta:
\(\dot{\mathbf q}=G\mathbf p\) and
\(\mathbf p=\mathcal G\dot{\mathbf q}\).
Equivalently, these are the Legendre relations associated with the kinetic
energy induced by the Riemannian metric \(\mathcal G\). Then the Newton form
of the LCCon dissipative particle equation becomes
\begin{equation}\label{LCCon-dissipative-particle-Newton-Riemannian-expanded}
\begin{split}
\ddot{\mathbf q}
+
\boldsymbol{\Gamma}(\dot{\mathbf q},\dot{\mathbf q})
+
\frac{\gamma}{1-z\psi}\dot{\mathbf q}
=
-
G\nabla_{\mathbf q}V
+
\frac{1}{1-z\psi}
G\mathbf{u}.
\end{split}
\end{equation}
Here the covector \(\mathbf u\) collects the locally conformal
correction terms and is given by
\begin{equation}\label{LCCon-u-term}
\mathbf u
=
(\boldsymbol{\varphi}(\mathbf q,z)\times
\mathcal G(\mathbf q)\dot{\mathbf q})
\times\dot{\mathbf q}
+
\boldsymbol{\varphi}(\mathbf q,z)
\left(
\frac{1}{2}
\dot{\mathbf q}\cdot
\mathcal G(\mathbf q)\dot{\mathbf q}
+
V(\mathbf q)
\right).
\end{equation}
The term \((1-z\psi)^{-1}G\mathbf u\) is the Riemannian realization of the
locally conformal force in the LCCon dissipative particle model. The vector
\(\mathbf u\) is naturally written in covector form, since it is constructed
from the \(d\mathbf q\)-component \(\boldsymbol{\varphi}\) of the Lee one-form,
while the inverse metric \(G\) raises the index and converts it into an
acceleration term. Its first component,
\((\boldsymbol{\varphi}\times\mathcal G\dot{\mathbf q})\times\dot{\mathbf q}\),
is a skew, gyroscopic-type correction produced by the conformal twisting,
whereas its second component is an energy-weighted force in the Lee direction.
The factor \((1-z\psi)^{-1}\) gives the additional conformal rescaling due to
the \(dz\)-component of the Lee form. Hence \((1-z\psi)^{-1}G\mathbf u\)
measures the deviation of the LCCon dissipative particle from the ordinary
contact dissipative particle on the Riemannian configuration manifold.

The action variable \(z\) satisfies
\begin{equation}\label{LCCon-dissipative-particle-z-Riemannian}
\dot z
=
\frac{1}{1-z\psi}
\left(
L_H(\mathbf q,\dot{\mathbf q},z)
+
z\,\boldsymbol{\varphi}\cdot\dot{\mathbf q}
\right).
\end{equation}
Here \(L_H(\mathbf q,\dot{\mathbf q},z)\) is the Herglotz Lagrangian \cite{cannarsa2019herglotz,LeonLainz2,Guenther} obtained
as the Legendre transformation of the contact Hamiltonian
\(H\) in \eqref{LCCon-dissipative-particle-H-Riemannian}, namely
\begin{equation}\label{Herglotz-Lagrangian-LCCon}
L_H(\mathbf q,\dot{\mathbf q},z)
=
\mathbf p\cdot\dot{\mathbf q}
-
H(\mathbf q,\mathbf p,z)
=
\frac{1}{2}\dot{\mathbf q}\cdot\mathcal G(\mathbf q)\dot{\mathbf q}
-
V(\mathbf q)
-
\gamma z .
\end{equation}
Thus, the $z$-equation is the locally conformal deformation of the usual
Herglotz equation. In particular, the term \(z\,\boldsymbol{\varphi}\cdot
\dot{\mathbf q}\) is the correction coming from the \(d\mathbf q\)-component
of the Lee form, while the denominator \(1-z\psi\) is produced by its
\(dz\)-component.

If \(\psi=0\), then the conformal denominator disappears and
\eqref{LCCon-dissipative-particle-Newton-Riemannian-expanded} reduces to
\begin{equation}\label{LCCon-dissipative-particle-Newton-Riemannian-psi-zero}
\begin{split}
\ddot{\mathbf q}
+
\boldsymbol{\Gamma}(\dot{\mathbf q},\dot{\mathbf q}) +\gamma\dot{\mathbf q}
=
-
G\nabla_{\mathbf q}V 
+
G\mathbf{u}.
\end{split}
\end{equation}
In this case, the action equation becomes
\(\dot z=L_H(\mathbf q,\dot{\mathbf q},z)
+z\,\boldsymbol{\varphi}\cdot\dot{\mathbf q}\).

Finally, when the Lee form vanishes, that is
\(\boldsymbol{\varphi}=0\) and \(\psi=0\), one has \(\mathbf u=0\), and the
ordinary contact dissipative particle equation on the Riemannian configuration
manifold is recovered:
\begin{equation}\label{contact-dissipative-particle-Riemannian}
\ddot{\mathbf q}
+
\boldsymbol{\Gamma}(\dot{\mathbf q},\dot{\mathbf q})
+
\gamma\dot{\mathbf q}
=
-
G\nabla_{\mathbf q}V.
\end{equation}
In the same limit, the action equation reduces to the standard Herglotz
relation \(\dot z=L_H(\mathbf q,\dot{\mathbf q},z)\).

\section{Summary and Discussion}\label{sec:summary-discussion}

We have developed a common local-to-global conformal framework for
cosymplectic and coorientable contact geometries and their Hamiltonian
dynamics. Non-degenerate odd pairs provide the shared algebraic setting,
while the Lichnerowicz--de~Rham differential encodes the conformal gluing
of local data. Propositions~\ref{prop:almost-poisson-Ham-closure}
and~\ref{prop:almost-Jacobi-Hamiltonian-obstruction} identify the tensors
governing the Hamiltonian bracket defects; the ordinary cosymplectic and
contact conditions eliminate these defects and recover Poisson and Jacobi
structures, respectively.

For cosymplectic geometry, we introduced an equal-weight LCCos formulation
and established its local-to-global characterization in
Proposition~\ref{prop:LCCos-local-to-global}. The associated bivector and
Lee vector field define a global Jacobi structure
(Proposition~\ref{prop:LCCos-Jacobi-structure}). With compatible Hamiltonians
carrying the same conformal weight as the two-form, the local Hamiltonian
vector fields glue directly without rescaling, and the resulting dynamics
agrees with direct Lichnerowicz--de~Rham deformation
(Proposition~\ref{prop:LCCos-Hamiltonian-equivalence}). These fields are
closed under the Jacobi--Lie bracket and generate an integrable
characteristic distribution
(Proposition~\ref{prop:LCCos-Hamiltonian-integrability}).

The conformally normalized LCCos evolution dynamics requires an additional
compatibility condition, expressed through a nowhere-vanishing
\(d_\varphi\)-closed scalar conformal unit \(\mathbb I\). Such a unit exists
precisely when the Lee form is exact. In this sector,
Proposition~\ref{prop:LCCos-evolution-equivalence} establishes the
equivalence of direct deformation and globalization of locally normalized
evolution fields, with
\(\mathcal E_H=X_H+\mathbb I\mathcal R^\tau\).
The time-dependent mass model in
Subsection~\ref{sec:LCCos-mass-example} illustrates the distinction between
Hamiltonian and evolution dynamics and the corresponding Lee-form
corrections.

For contact geometry, the twisted odd pair
\(\Omega=d_\varphi\eta\) yields a global Jacobi structure after correcting
the odd-pair Reeb vector field
(Proposition~\ref{prop:LCCon-global-Jacobi}).
Proposition~\ref{prop:LCCon-underlying-contact} identifies this structure
with the canonical Jacobi structure of the ordinary contact form \(\eta\):
the odd-pair bivector is the contact bivector, and the corrected vector
field is the contact Reeb field.
Proposition~\ref{prop:LCCon-Hamiltonian-equivalence} then establishes that
direct deformation and conformal globalization of compatible local contact
Hamiltonian systems yield the same global vector field. The dissipative
particle model in Subsection~\ref{sec:LCCon-particle-example} makes explicit
the resulting coordinate corrections, including gyroscopic-type,
energy-dependent, and velocity-linear contributions.

The following directions will be investigated in future studies.

\noindent
\textbf{Kinetic Theories for LCCos and LCCon Dynamics.}
Particle dynamics can be lifted to kinetic equations through suitable
geometric lifts and vertical representatives
\cite{EsGu11,EsGu12,EsGrMiGu19}. Related kinetic formulations have already
been developed for contact geometry \cite{esen2025conformal} and locally
conformally symplectic dynamics \cite{esen2023conformal}. A natural
continuation of the present work is to construct kinetic models for both the
LCCos and LCCon dynamics developed here. Such formulations would incorporate
explicit time dependence, Lee-form corrections, and contact-type dissipation
into the geometric description of ensemble dynamics.

\noindent
\textbf{Geometric Hamilton--Jacobi Theory for LCCos and LCCon Dynamics.}
Geometric Hamilton--Jacobi theory
\cite{carinena2006geometric,esen2022reviewing,roman2021overview} has been
extended to locally conformally symplectic mechanics
\cite{EsenLeonSarZaj1}, different-weight LCCos systems
\cite{AtEsLe23}, \(k\)-cosymplectic field theories
\cite{HJ-co-k-field}, locally conformally \(k\)-symplectic geometry
\cite{EsLeSaZa-k-sympl}, and multisymplectic geometry
\cite{EsLeSaZa-Cauchy}. The equal-weight formulation introduced here calls
for a corresponding reformulation of the LCCos Hamilton--Jacobi theory and,
subsequently, its extension to equal-weight locally conformally
\(k\)-cosymplectic and LCCon dynamics.

\noindent
\textbf{Locally Conformally Cocontact Geometry and Dynamics.}
Cocontact geometry combines explicit time dependence with contact-type
dissipation
\cite{azuaje2023canonical,LeGaGrMuRi23,Leon-cocontact23,gaset2023symmetries,MaNo17,
RiTo23}. A separate future study will combine the LCCos and
LCCon structures developed here into an LCCoco framework based on a closed
Lee form and the Lichnerowicz--de~Rham differential. It will address the
globalization of local cocontact structures, the transformation of the two
Reeb vector fields and the cocontact bivector, the resulting Jacobi
structure, and the associated Hamiltonian and evolution dynamics. The
time-dependent mass and dissipative particle models considered in this paper
will provide a natural mechanical application of this unified framework.

\section*{AI Use Statement}

The authors used GPT (OpenAI) as an auxiliary tool for language editing during manuscript preparation and for cross-checking selected algebraic computations. The authors reviewed and approved all content in its final form and take full responsibility for the manuscript.

\section*{Acknowledgment}
 B. Ateşli and A. Çatal-Özer gratefully acknowledge the support of TÜBİTAK through the 2218 Postdoctoral Research Fellowship Program, Project No. 124C484.

 \appendix

 \section{Appendix}

 \subsection{Schouten--Nijenhuis Algebra} \label{sec:SN}

For multivector fields $\mathcal{X}$ of order $m \geq 1$ and   $\mathcal{Y}$ of order $n \geq 1$ given by
\begin{equation} \label{decompi}
\mathcal{X}=X_1\wedge \dots\wedge X_m, \qquad \mathcal{Y}=Y_1\wedge \dots\wedge Y_n,
\end{equation}
the Schouten--Nijenhuis bracket in terms of the Jacobi-Lie bracket is as follows
\begin{equation} \label{SN-def1}
[\mathcal{X},\mathcal{Y}]_{SN}:=\sum_{i,j}(-1)^{i+j}[X_i,Y_j]\wedge X_1\wedge \dots\wedge\widehat X_{i} \wedge \dots \wedge X_m\wedge Y_1\wedge \dots \wedge \widehat Y_{j} \wedge \dots \wedge Y_n, 
\end{equation}
where $\widehat{X}_{i}$ and $\widehat{Y}_{j}$ denote the omission of the elements $X_i$ and $Y_j$, respectively, from the sequences. Note that the resulting multivector field $[\mathcal{X},\mathcal{Y}]_{SN}$ is of order $m+n-1$. For the case $m=1$ and $n=1$, one has two vector fields and the Schouten--Nijenhuis bracket becomes the Jacobi-Lie bracket: \begin{equation}
[X,Y]_{SN}=[X,Y]=\mathcal{L}_X Y.
\end{equation}

For the case $m \geq 1$ and $n=0$, that is, if there is a multivector field $\mathcal{X}$ of order $m\geq 1$ and a function $F$, the Schouten--Nijenhuis bracket is defined as 
\begin{equation}\label{X-F-SN}
    [\mathcal{X},F]_{SN} := \sum^m_{i=1}(-1)^{i+1}X_i(F) X_1\wedge \dots\wedge\widehat X_{i} \wedge \dots \wedge X_m. 
\end{equation}
Note that $[\mathcal{X},F]_{SN}$ is a multivector field of order $m-1$.
If $m=1$, that is, when $\mathcal{X}$ is a vector field, then
\begin{equation}
 [X,F]_{SN}=[X,F]=X(F)=\mathcal{L}_X F,
\end{equation}
where $X(F)$ is the directional derivative, that is, the Lie derivative of $F$ in the direction of $X$.

For $m=n=0$, one has two functions, say $F$ and $H$. The Schouten--Nijenhuis bracket is then defined to be identically zero, that is
\begin{equation}
    [F,H]_{SN}:=0.
\end{equation}

 Consider a multivector field $\mathcal{X}$ of order $m \geq 1$ and a multivector field $\mathcal{Y} $ of order $n\geq 1$. The Schouten--Nijenhuis bracket is a supercommutative bracket satisfying
\begin{equation}\label{SN-commute}
    [\mathcal{Y},\mathcal{X}]_{SN} = -(-1)^{(m-1)(n-1)} [\mathcal{X},\mathcal{Y}]_{SN}.
\end{equation}
Note that two multivector fields of even order commute, while any other combination of orders gives rise to anticommutation of such multivectors. Now, consider a multivector field $\mathcal{Z}$ of positive order. Then, the Schouten--Nijenhuis bracket satisfies the following identity:
\begin{equation}\label{SN-id}
    [\mathcal{X},\mathcal{Y} \wedge \mathcal{Z}]_{SN} = [\mathcal{X},\mathcal{Y}]_{SN} \wedge \mathcal{Z} + (-1)^{(m+1)n}\mathcal{Y} \wedge [\mathcal{X},\mathcal{Z}]_{SN}.
\end{equation}
For a function $F$, an identity similar to \eqref{SN-id} can be written as
\begin{equation}\label{X-FY}
     [\mathcal{X},F\mathcal{Y}]_{SN} = F[\mathcal{X},\mathcal{Y}]_{SN} + (-1)^{m+1} [\mathcal{X},F] _{SN}\wedge \mathcal{Y},
\end{equation}
where $\mathcal{X}$ is a multivector field of order $m$.

\subsection{Right Contraction Mapping (Interior Derivative)} \label{sec:intder} 

The (right) contraction of the multivector field $\mathcal{X}$ in \eqref{decompi}  with a one-form $\theta$ is defined as
\begin{equation}\label{right-cont}
    \iota_{\theta}\mathcal{X} := \sum_{i=1}^{m} (-1)^{i+m} \langle \theta, X_i \rangle X_1 \wedge \dots \wedge \widehat{X_i} \wedge \dots \wedge X_m.
\end{equation}
By taking $\theta = dF$ as an exact one-form, we can express the Schouten--Nijenhuis bracket of a multivector field $\mathcal{X}$ with a smooth function $F$ as in \eqref{X-F-SN}:
\begin{equation} 
  [\mathcal{X}, F]_{SN} = (-1)^{m+1} \iota_{dF} \mathcal{X}.
\end{equation}
Using the right contraction mapping, the identity exhibited in \eqref{X-FY} becomes
\begin{equation}\label{X-FY-2}
    [\mathcal{X}, F \mathcal{Y}]_{SN} = F [\mathcal{X}, \mathcal{Y}]_{SN} + \iota_{dF} \mathcal{X} \wedge \mathcal{Y}.
\end{equation}

We now show how the right contraction mapping distributes over the wedge product. We consider once more the multivectors $\mathcal{X}$ and $\mathcal{Y}$ given in \eqref{decompi}. We have 
\begin{equation}\label{iota-SN-m-n}
     \iota_{dF} (\mathcal{X} \wedge \mathcal{Y}) = (-1)^n\iota_{dF}\mathcal{X} \wedge \mathcal{Y} + \mathcal{X} \wedge \iota_{dF} \mathcal{Y} .
\end{equation}
For $m=1$ and $n=2$,
\begin{equation}\label{iota-SN-2-2-X}
 \iota_{dF} (X \wedge \mathcal{Y}) =X(F) \mathcal{Y} + X\wedge  \iota_{dF} \mathcal{Y}.
\end{equation}
For $m=n=2$,
\begin{equation}\label{iota-SN-2-2-X-}
 \iota_{dF} (\mathcal{X} \wedge \mathcal{Y}) = \iota_{dF}\mathcal{X} \wedge \mathcal{Y} + \mathcal{X} \wedge \iota_{dF} \mathcal{Y} .
\end{equation}
For $\alpha\in\Gamma^k(\mathcal{Q})$, $\beta\in\Gamma^l(\mathcal{Q})$, we have
\begin{equation}
    \iota_{\alpha \wedge \beta} X = \iota_\alpha\iota_\beta X.
\end{equation}
In particular, for $F_1,...,F_k$ in $C^\infty(\mathcal{Q})$,
\begin{equation}
     \iota_{dF_1 \wedge \dots  \wedge dF_k} X = \iota_{dF_1}\dots \iota_{dF_k} X .
\end{equation}

For a vector field $X$ and a bivector field $\mathcal{Y}$, the right contraction $\iota_{dF}$ distributes over the Schouten--Nijenhuis bracket as follows:
\begin{equation}
    \iota_{dF}[X,\mathcal{Y}]_{SN} = [\mathcal{Y},\iota_{dF}X]_{SN} + [X,\iota_{dF}\mathcal{Y}]_{SN}.
\end{equation}
For two bivector fields $\mathcal{X} = X_1 \wedge X_2$ and $\mathcal{Y} = Y_1 \wedge Y_2$, the right contraction $\iota_{dF}$ distributes over the Schouten--Nijenhuis bracket as follows:
\begin{equation}\label{iota-SN-2-2+}
    \iota_{dF}([\mathcal{X},\mathcal{Y}]_{SN}) = [\mathcal{Y},\iota_{dF}(\mathcal{X})]_{SN} + [\mathcal{X},\iota_{dF}(\mathcal{Y})]_{SN},
\end{equation}
for a function $F$. Note that, if $\mathcal{X}=\mathcal{Y}$ is of order $2$, then
\begin{equation}\label{iota-SN-2-2-X+}
 \iota_{dF}([\mathcal{X},\mathcal{X}]_{SN}) = 2[\mathcal{X},\iota_{dF}(\mathcal{X})]_{SN}. 
\end{equation}
In general, for $\mathcal{X}= X_1 \wedge \dots X_m$ and $\mathcal{Y}= Y_1 \wedge \dots Y_n$, the right contraction $\iota_{dF}$ distributes over the Schouten--Nijenhuis bracket as
\begin{equation}
    \iota_{dF}([\mathcal{X},\mathcal{Y}]_{SN}) = (-1)^n [\mathcal{Y},\iota_{dF}\mathcal{X}]_{SN} + [\mathcal{X},\iota_{dF}\mathcal{Y}]_{SN}.
\end{equation}

 \bibliographystyle{abbrv}
\bibliography{references}

\end{document}